\documentclass[letterpaper,12pt]{amsart}

\usepackage[margin=1.14in]{geometry}
\usepackage{amsmath,amsthm,amssymb,mathtools,amsfonts,mathrsfs}

\usepackage{tikz-cd}
\usepackage[all]{xy}
\usepackage{array,longtable,makecell,diagbox,booktabs,tabularx}
\usepackage{setspace}
\usepackage[
  breaklinks,
  colorlinks=true,
  citecolor=teal,
  linkcolor=teal,
  urlcolor=teal,
  pagebackref,
  hyperindex
]{hyperref}
\usepackage[alphabetic]{amsrefs}

\numberwithin{equation}{section}

\theoremstyle{plain}
\newtheorem{thm}{Theorem}[section]
\newtheorem{prop}[thm]{Proposition}
\newtheorem{lem}[thm]{Lemma}
\newtheorem{cor}[thm]{Corollary}
\newtheorem{conj}[thm]{Conjecture}
\newtheorem{ques}[thm]{Question}

\theoremstyle{definition}
\newtheorem{defn}[thm]{Definition}
\newtheorem{rmk}[thm]{Remark}
\newtheorem{example}[thm]{Example}

\newcommand{\C}{\mathbb C}

\newcommand{\PP}{\mathbb P}

\newlength{\notationcolwidth}
\newlength{\potentialcolwidth}

\newcommand{\Aut}{\operatorname{Aut}}
\newcommand{\ev}{\operatorname{ev}}

\newcommand{\vir}{\operatorname{vir}}

\newcommand{\Res}{\mathrm{Res}}

\newcommand{\Vir}{\operatorname{Vir}}
\newcommand{\id}{\operatorname{id}}

\title{Virasoro operators for relative Gromov--Witten theory}
\author{Shuai Guo}
\address{School of Mathematical Sciences, Peking University, Beijing 100871, China}
\email{guoshuai@math.pku.edu.cn}

\author{Longting Wu}
\address{Department of Mathematics \& SUSTech International Center for Mathematics, Southern University of Science and Technology,
Shenzhen, Guangdong 518055, China}
\email{wult@sustech.edu.cn}

\author{Lei You}
\address{School of Mathematical Sciences, Peking University, Beijing 100871, China}
\email{youlei26@stu.pku.edu.cn}
\date{}

\begin{document}
\raggedbottom

\begin{abstract}
Let $D$ be a smooth divisor in a smooth projective variety $X$.
We construct relative Virasoro operators acting on the relative
descendant potential of $(X,D)$. Each operator is the
logarithmic Virasoro operator of $(X,D)$ plus a boundary correction
determined only by $D$ and its normal bundle in $X$.
Assuming the Virasoro conjecture for $D$, we prove that the
relative constraints are equivalent to the absolute constraints
for $X$. Under the same assumption,
we establish compatibility with the degeneration formula for
two-component projective semistable degenerations satisfying
a completeness condition. Consequently, the Virasoro conjecture
for the two components and the gluing divisor implies the
conjecture for the general fiber. As an application, we prove
the Virasoro conjecture for Enriques surfaces in every integral
curve class.
\end{abstract}

\maketitle

\tableofcontents

\section{Introduction}

\subsection{Absolute Gromov-Witten theory and Virasoro operators}

Let $X$ be a smooth projective variety over $\C$ of complex dimension
$d$. We use 
$\overline{\mathcal M}_{g,n}(X,\beta)$ to denote the moduli space of genus $g$, $n$-pointed
stable maps to $X$ with curve class $\beta\in H_2(X,\mathbb Z)$.

For nonnegative integers $m_1,\ldots,m_n$ and cohomology classes
$\gamma_1,\ldots,\gamma_n\in H^*(X,\mathbb Q)$, the corresponding
descendant Gromov--Witten invariant of $X$ is defined by
\begin{equation}\label{eq:absolute-gw}
 \left\langle\prod_{i=1}^n\gamma_i\psi_i^{m_i}\right\rangle^X_{g,\beta}
 :=\int_{[\overline{\mathcal M}_{g,n}(X,\beta)]^{\vir}}
       \prod_{i=1}^n\psi_i^{m_i}\ev_i^*\gamma_i
\end{equation}
where $\ev_i$ and $\psi_i$ denote the evaluation map and the psi class
at the $i$th marked point.

We choose a Hodge-homogeneous basis
$\{\phi_a\}$ of $H^*(X,\mathbb C)$, indexed by a
finite set $A$ with distinguished element $0$, such that
$\phi_0=1_X\in H^0(X,\mathbb C)$ and
$\phi_a\in H^{p_a,q_a}(X)$.
Let $Q^\beta$ denote the Novikov variable associated to the curve class
$\beta$, and let $t_m^a$ be the super formal variable corresponding to
the descendant insertion $\phi_a\psi^m$. We write
\[
  \mathbf t=\mathbf t(\psi):=\sum_{m\geq0}\sum_{a}t_m^a\psi^m\phi_a.
\]
Using the convention
\[
  \left\langle e^{\mathbf t}\right\rangle^X_{g,\beta}
  :=
  \sum_{n\geq0}\frac{1}{n!}
  \left\langle\prod_{i=1}^n\mathbf t(\psi_i)\right\rangle^X_{g,\beta},
\]
we define the genus-$g$ connected descendant potential by
\begin{equation}\label{eq:absolute-connected-potential}
F_g^X(\mathbf t,Q)
  =
  \sum_\beta Q^\beta
  \left\langle e^{\mathbf t}\right\rangle^X_{g,\beta}.
\end{equation}
We further define the total connected descendant potential and the disconnected
descendant potential by
\begin{equation}\label{eq:absolute-total-potential}
 F^X(\mathbf t,Q,\hbar)=\sum_{g\ge0}\hbar^{g-1}F_g^X(\mathbf t,Q),\qquad
 Z^X(\mathbf t,Q,\hbar)=\exp\bigl(F^X(\mathbf t,Q,\hbar)\bigr).
\end{equation}

We next recall the absolute Virasoro operators acting on the descendant
partition function $Z^X$. These operators were introduced by
Eguchi, Hori, and Xiong \cite{EHX}. Their extension to arbitrary compact K\"ahler targets, due to S.~Katz, appears in \cite{EJX} and is reviewed
in \cite{G99}.
Set
\[
 \eta^{X}_{ab}=\int_X\phi_a\cup\phi_b,\qquad
 b_a=p_a+\frac{1-d}{2},\qquad
 \rho(\gamma)=c_1(X)\cup\gamma.
\]
Write $\eta_X^{ab}$ for the entries of the inverse matrix of
$(\eta^{X}_{ab})$, and set
\begin{align*}
 \rho^i\phi_a&=\sum_b(\rho^i)_a{}^b\phi_b,\\
 (\rho^i)_{ab}&=\int_X\phi_a\cup\rho^i\phi_b,\qquad
 (\rho^i)^{ab}=\sum_c\eta_X^{ac}(\rho^i)_c{}^b.
\end{align*}

Define the dilaton-shifted coefficients
\[
 \widetilde t_m^a=t_m^a-\delta_{m,1}\delta_{a,0},\qquad
\partial_{m,a}
  :=\frac{\partial}{\partial t_m^a}
  \quad (m\geq0)
\]
and set $\partial_{m,a}=0$ for $m<0$.
For $k\ge0$ and $0\le i\le k+1$, let
\begin{equation}\label{eq:virasoro-coefficients}
 [x]^k_i=e_{k+1-i}(x,x+1,\ldots,x+k)
\end{equation}
where $e_r$ is the elementary symmetric polynomial and $e_0=1$.

For every $k\ge -1$, the absolute Virasoro operator
is defined by
\begin{align}
 L_k^X={}&\sum_{m\ge0}\sum_{i=0}^{k+1}\sum_{a,b}
 [m+b_a]^k_i(\rho^i)_a{}^b
 \widetilde t_m^a\partial_{m+k-i,b}\notag\\
 &+\frac{\hbar}{2}\sum_{i=0}^{k-1}\sum_{m=0}^{k-i-1}
 \sum_{a,b}(-1)^{m+1}[-m-b_a]^k_i
 (\rho^i)^{ab}\partial_{m,a}\partial_{k-i-1-m,b}
 \label{eq:absolute-virasoro}\\
 &+\frac1{2\hbar}\sum_{a,b}(\rho^{k+1})_{ab}t_0^bt_0^a
 +\frac{\delta_{k,0}}{48}
 \int_X\left((3-d)c_d(X)
 -2c_1(X)c_{d-1}(X)\right).\notag
\end{align}
Here we set
$[x]^{-1}_0=1$, $[x]^{-1}_i=0$ for $i>0$, and interpret empty sums as zero. These operators satisfy the Virasoro bracket relations:
\begin{equation}\label{eq:absolute-virasoro-commutator}
 [L_k^X,L_l^X]=(k-l)L_{k+l}^X.
\end{equation}
\begin{conj}[Virasoro conjecture]\label{conj:absolute-virasoro}
For every smooth projective variety $X$, the descendant partition
function satisfies
\begin{equation*}
 L_k^XZ^X=0\qquad(k\ge-1).
\end{equation*}
\end{conj}

In all genera, two important approaches to the Virasoro conjecture
are quantization and degeneration. The quantization formalism of
Givental \cite{Giv5}, together with the classification of semisimple
cohomological field theories by Teleman \cite{T12}, establishes the
constraints for targets with generically semisimple quantum cohomology.
Using the quantization formalism, Coates, Givental, and Tseng
\cite{CGT} further proved that the Virasoro constraints hold for
the total space of a toric bundle if and only if they hold for its base.

A different approach uses degeneration and relative Gromov--Witten
theory. Okounkov and Pandharipande \cite{OP3} proved the constraints
for all nonsingular projective curves, including curves of positive
genus, whose quantum cohomology is not semisimple.
Hu, Jin, Ke, Lu, Tseng, and Wu \cite{HJKLTW} extended this approach
to smooth projective effective orbifold curves.
In these works, relative Virasoro constraints and their compatibility
with degeneration play a central role in reducing the problem to
simpler relative geometries.
OpenAI \cite{OpenAIci} recently announced a degeneration-based proof of the
Virasoro constraints for smooth projective complete intersections,
without constructing relative Virasoro operators.

\subsection{Relative Gromov-Witten theory}
Let $D\subset X$ be a smooth divisor, assumed for now to be connected.
Fix a Hodge-homogeneous basis $\{\delta_b\}$ of $H^*(D,\mathbb C)$
with $\delta_0=1_D\in H^0(D,\mathbb C)$.
A relative condition is a
cohomology-weighted partition
\[
 \lambda=\bigl((\mu_1,\delta_{b_1}),\ldots,(\mu_\ell,\delta_{b_\ell})\bigr),
 \qquad \mu_j\in\mathbb Z_{>0},
\]
where $\mu_j$ is the contact order and $\delta_{b_j}$ is the
chosen basis element at the $j$th relative marking.
We denote
\[
|\lambda|=\sum_{j=1}^{\ell}\mu_j,\qquad  \ell(\lambda)=\ell.
\]

Let $\overline{\mathcal M}^{\bullet}_{g,n}(X/D,\beta,\mu)$ denote the moduli
space of possibly disconnected relative stable maps with $n$ ordinary markings
and $\ell$ labelled relative markings of contact orders
$\mu=(\mu_1,\ldots,\mu_\ell)$.
The maps take values in expanded
targets of the pair $(X,D)$, and their curve class after contraction
to $X$ is $\beta$ \cites{Jun1,Jun2}. Compatibility with the curve class $\beta$ requires
$|\lambda|=[D]\cdot\beta$. 
Here $\bullet$ allows disconnected domains, and
$g=1-\chi(\mathcal O_C)$ is the arithmetic genus of the domain $C$.
Thus, if $C$ has $c$ connected components of genera $g_1,\ldots,g_c$,
then $g-1=\sum_{i=1}^c(g_i-1)$; in particular, $g$ may be negative.
Stable degree-zero connected components are allowed.

This moduli space carries a virtual
fundamental class of complex dimension
\[
 (1-g)(d-3)+\int_\beta c_1(X)+n+\ell(\lambda)-|\lambda|.
\]
The ordinary and relative evaluation maps take values in $X$ and
$D$, respectively, and are denoted by $\ev_i$ and $\ev_j^{\mathrm{rel}}$.

For ordinary descendant insertions, define the disconnected relative
Gromov--Witten invariants by
\begin{equation}\label{eq:relative-gw}
 \left\langle\prod_{i=1}^n\gamma_i\psi_i^{m_i}
       \mathrel{\big|}\lambda\right\rangle^{\bullet,X/D}_{g,\beta}
 :=\int_{[\overline{\mathcal M}^{\bullet}_{g,n}(X/D,\beta,\mu)]^{\vir}}
 \prod_{i=1}^n\psi_i^{m_i}\ev_i^*(\gamma_i)
 \prod_{j=1}^{\ell}(\ev_j^{\mathrm{rel}})^*(\delta_{b_j}).
\end{equation}
Here $\psi_i$ are the cotangent line classes at the ordinary
markings. The relative markings are labelled, so the bracket is defined
without dividing by $|\Aut(\lambda)|$.
Using the same formal insertion $\mathbf t=\mathbf t(\psi)$ as in the absolute theory, we define
the disconnected relative descendant potential by
\begin{equation}\label{eq:relative-disconnected-potential}
 Z^{X/D}(\mathbf t,Q,\hbar\mid\lambda)
 :=\sum_{\beta:\,[D]\cdot\beta=|\lambda|}\sum_{g\in\mathbb Z}
 Q^\beta\hbar^{g-1}
 \left\langle e^{\mathbf t}\mathrel{\big|}\lambda\right\rangle^{\bullet,X/D}_{g,\beta},
\end{equation}
where the exponential denotes the sum of insertions
$\prod_{i=1}^n\mathbf t(\psi_i)$ with weight $1/n!$,
while the relative condition $\lambda$ is fixed.
The empty domain contributes $1$ when $\beta=0$ and
$\lambda=\varnothing$, with no ordinary markings.

If $D=\coprod_{s=1}^r D_s$ is disconnected, the same definitions apply
after specifying the component $D_{s(j)}$ containing the image of each
relative marking. We choose the basis of $H^*(D,\mathbb C)$ to be the union of homogeneous
bases of $H^*(D_s,\mathbb C)$ so that $\delta_{b_j}$ is supported on
$D_{s(j)}$. The compatibility condition becomes
\[
 \sum_{j:\,s(j)=s}\mu_j=[D_s]\cdot\beta,\qquad 1\le s\le r,
\]
and the sum in \eqref{eq:relative-disconnected-potential} is restricted
to curve classes satisfying all these conditions.

\subsection{Relative Virasoro operators}

For a nonsingular projective curve $X$ and a reduced divisor
$D=p_1+\cdots+p_m$, Okounkov and Pandharipande \cite{OP3}
constructed relative Virasoro operators by replacing the tangent bundle $T_X$
in the absolute construction by the logarithmic tangent bundle
$T_X(-\log D)$. In terms of the coefficients, this amounts to
the replacement
\[
 c_1(X)\longmapsto c_1\bigl(T_X(-\log D)\bigr)=c_1(X)-[D].
\]
Their operators $L_k^{X/D}\, (k\geq -1)$ act on the ordinary descendant variables
and are independent of relative conditions. They still satisfy the Virasoro bracket relation:
\[
 [L_k^{X/D},L_l^{X/D}]=(k-l)L_{k+l}^{X/D}.
\]
They annihilate the relative potential
\cite{OP3}*{Theorem~3}:
\[
 L_k^{X/D}Z^{X/D}=0,
 \qquad k\ge-1.
\]
These constraints are compatible with the degeneration formula:
gluing the relative theories along matching ramification profiles
preserves the Virasoro constraints.

This leads naturally to the following question in relative Gromov-Witten theory:
\begin{ques}
For a smooth projective variety $X$ of arbitrary dimension and a smooth
divisor $D\subset X$, can one construct relative Virasoro operators
extending those of Okounkov-Pandharipande and compatible with the
degeneration formula?
\end{ques}

Our aim is to answer this question. In higher dimensions, simply
replacing $T_X$ with $T_X(-\log D)$ in the absolute construction
is generally insufficient. Additional boundary corrections may
mix different relative conditions, leading us to assemble the
relative potentials into a single vector-valued generating function.

Assume first that $D$ is connected, and use the chosen homogeneous
basis $\{\delta_b\}$ of $H^*(D,\mathbb C)$.
For each $\mu\ge1$, introduce a generator $p_{\mu,b}$
recording contact order $\mu$ and basis index $b$,
and set
\[
\mathcal H_D
=\operatorname{Sym}_{\mathrm{super}}\!\left(
  \bigoplus_{\mu\ge1,b}\mathbb C\,p_{\mu,b}
\right).
\]
Each generator $p_{\mu,b}$ has contact degree $\mu$ and parity
$|\delta_b|:=\deg_{\mathbb R}(\delta_b)\bmod 2$.

Let $(\mathcal H_D)_d$ denote the subspace of total contact degree $d$,
and define
\[
 \widehat{\mathcal H}_D=\prod_{d\ge0}(\mathcal H_D)_d.
\]
We call $\mathcal H_D$ the boundary state space and
$\widehat{\mathcal H}_D$ its completion with respect to total contact degree.

Fix an ordering of the generators. A relative condition
$\lambda=((\mu_1,\delta_{b_1}),\ldots,(\mu_\ell,\delta_{b_\ell}))$ determines
the ordered monomial
\[
 p_\lambda=\prod_{j=1}^{\ell}p_{\mu_j,b_j},
 \qquad p_{\varnothing}=1.
\]
The monomials with no repeated odd generator form a basis of
$\mathcal H_D$.
We use the following exponential notation for the sum over
relative conditions:
\[
 \left\langle -\mathrel{\big|}e^{\mathbf p}\right\rangle^{\bullet,X/D}_{g,\beta}
 :=\sum_{\lambda:\,|\lambda|=[D]\cdot\beta}
 \frac{1}{|\Aut(\lambda)|}
 \left\langle -\mathrel{\big|}\lambda\right\rangle^{\bullet,X/D}_{g,\beta}p_\lambda.
\]
All products follow the Koszul sign convention: interchanging
homogeneous factors of parities $\epsilon$ and $\epsilon'$
introduces the sign $(-1)^{\epsilon\epsilon'}$.
The disconnected relative potential is then
\begin{equation}\label{eq:relative-total-potential}
\begin{aligned}
 \mathbf Z^{X/D}(\mathbf t,Q,\hbar;\mathbf p)
 &:=\sum_{\beta,\,g\in\mathbb Z}Q^\beta\hbar^{g-1}
 \left\langle e^{\mathbf t}\mathrel{\big|}e^{\mathbf p}
 \right\rangle^{\bullet,X/D}_{g,\beta}\\
 &=\sum_\lambda  \frac{1}{|\Aut(\lambda)|}Z^{X/D}(\mathbf t,Q,\hbar\mid\lambda)p_\lambda.
\end{aligned}
\end{equation}
The potential $\mathbf Z^{X/D}$ may be viewed either as a
generating function in the boundary variables or as a
$\widehat{\mathcal H}_D$-valued formal function.

We seek relative Virasoro operators acting on this potential.
In the basis $\{p_\lambda\}$, they are operator-valued matrices
that may mix different relative conditions. Compatibility with
the degeneration formula is expressed through the gluing pairing
of these potentials.
\begin{rmk}
If $D=\coprod_{s=1}^r D_s$ is disconnected, the same definition of
$\mathcal H_D$ applies, using
$H^*(D,\mathbb C)=\bigoplus_s H^*(D_s,\mathbb C)$.
There is a natural isomorphism
\[
 \mathcal H_D\cong\bigotimes_{s=1}^r\mathcal H_{D_s},
\]
where the tensor product is taken in the category of superalgebras.
The state space is naturally multigraded by the contact degrees
$\mathbf d=(d_1,\ldots,d_r)$ along the connected components:
a generator $p_{\mu,b}$ with $\delta_b \in H^*(D_s,\mathbb C)$ has
degree $\mu$ in the $s$th component and zero in the others. We denote
\[
 \widehat{\mathcal H}_D
 =\prod_{\mathbf d\in\mathbb Z_{\ge0}^r}(\mathcal H_D)_{\mathbf d}.
\]
This agrees with completion by total contact degree
$d=d_1+\cdots+d_r$. The coefficient of $Q^\beta$ in $\mathbf Z^{X/D}$
lies in multidegree $([D_1]\cdot\beta,\ldots,[D_r]\cdot\beta)$.
\end{rmk}

We can now state our main result. Let $\operatorname{Vir}(X)$ denote the assertion that
\[
L_k^X Z^X=0, \qquad \text{for all } k\ge -1.
\]
To define \(L_k^{X/D,\log}\), we use \eqref{eq:absolute-virasoro} with the same Hodge grading and dilaton shift, replacing \(T_X\) by \(T_X^{\log}:=T_X(-\log D)\) in every Chern class. Its Chern operator is therefore
\[
 c_1(T_X^{\log})\cup=(c_1(X)-[D])\cup,
\]
and its scalar term is $\delta_{k,0} \varepsilon^{X/D,\log}$, where
\begin{equation}\label{eq:logarithmic-scalar}
 \varepsilon^{X/D,\log}
 :=\frac1{48}\int_X
   \bigl((3-d)c_d(T_X^{\log})-2c_1(T_X^{\log})c_{d-1}(T_X^{\log})\bigr),
 \qquad d=\dim_{\mathbb C}X.
\end{equation}
On $\widehat{\mathcal H}_D$-valued functions, this operator acts componentwise.
We write $\varepsilon^X$ for the scalar coefficient of $L_0^X$
in \eqref{eq:absolute-virasoro}.

Since replacing $c_1(X)\cup$ by $(c_1(X)-[D])\cup$ does not
change the scalar contribution in the commutator calculation
of \cite{G99}*{Theorem~2.2},
for $k,l\ge-1$ and $k+l\ge-1$,
\begin{equation*}
 [L_k^{X/D,\log},L_l^{X/D,\log}]
 =(k-l)L_{k+l}^{X/D,\log}
 -(k-l)\delta_{k+l,0}
 (\varepsilon^{X/D,\log}-\varepsilon^X)\operatorname{id}.
\end{equation*}
Only the pairs $(k,l)=(-1,1),(1,-1)$ can have this additional
scalar. For curves both constants vanish, in agreement with the
relative operators recalled above. We use the logarithmic scalar
because it is additive under degeneration, as shown in
\eqref{eq:logarithmic-scalar-additivity}.

\begin{thm}\label{thm:relative-virasoro-main}
Let $X$ be a smooth projective variety and
$D=\coprod_{s=1}^r D_s\subset X$ a smooth divisor with connected
components $D_s$. We use the same Hodge grading convention as in the
absolute theory.
Write $i_s:D_s\hookrightarrow X$ and $N_s=N_{D_s/X}$.
There exist relative Virasoro operators $L_k^{X/D}$ for all $k\ge-1$
with the following properties.
\begin{enumerate}
\item They have the form
\[
 L_k^{X/D}=L_k^{X/D,\log}
       +\sum_{s=1}^r\Delta L_k^{D_s,N_s}(i_s^*\mathbf t).
\]
Each correction $\Delta L_k^{D_s,N_s}$ is a formal boundary matrix
determined only by the local pair $(D_s,N_s)$ and acts on the
corresponding boundary factor. Its dependence on the ordinary
insertions is through $i_s^*\mathbf t$, with formal dependence on the genus
and curve variables understood.

\item If $\operatorname{Vir}(D_s)$ holds for every $s$, the operators
$L_k^{X/D}$ preserve $\mathscr S_{X/D}$, defined by
\eqref{eq:insertion-solution-space} and containing $\mathbf Z^{X/D}$.
On $\mathscr S_{X/D}$, they satisfy
\[
 [L_k^{X/D},L_l^{X/D}]=(k-l)L_{k+l}^{X/D},
 \qquad(k,l\ge-1,\ k+l\ge-1).
\]

\item Under the same assumption, the absolute and relative Virasoro
constraints are equivalent:
\[
 \operatorname{Vir}(X)\quad\Longleftrightarrow\quad
 \operatorname{Vir}(X/D)
\]
where $\operatorname{Vir}(X/D)$ denotes
$L_k^{X/D}\mathbf Z^{X/D}=0$ for every $k\ge-1$.
\end{enumerate}
\end{thm}

\begin{rmk}
The first two relative operators are
\[
 L_{-1}^{X/D}=L_{-1}^{X/D,\log},\qquad
 L_0^{X/D}=L_0^{X/D,\log}+\mathfrak b.
\]
Here $\mathfrak b$ is a diagonal boundary grading operator; its formula
is given in Section~\ref{subsec:relative-string-dilaton}.
In Section~\ref{sec:examples}, we illustrate the corrections: for a
curve $X$, they all vanish and
$L_k^{X/D}=L_k^{X/D,\log}$ for every $k\ge-1$
(Section~\ref{subsec:curve-example}), while in higher dimensions they
can be nonzero. The commutation relations in part~(2) may not hold on
the full space of relative state-valued series; see
Section~\ref{subsec:p2-virasoro}.
\end{rmk}

The construction in Section~\ref{sec:relative-virasoro-construction}
uses deformation to the normal cone $X\rightsquigarrow X\cup_D Y$,
where $Y=\PP_D(N\oplus\mathcal O_D)$ and $N=N_{D/X}$.
The pair $(Y,D_0)$, with $D_0$ the zero section, is the cap.
Comparing the absolute Virasoro action with the logarithmic actions
on the two components yields a correction $\mathcal A_k^{D,N}$
expressed by rubber invariants, as in
\eqref{eq:relative-cap-defect-identity}.
Writing $\mathsf C$ for contraction with the cap potential and
$\mathsf R$ for the cap reconstruction operator of
Section~\ref{subsec:cap-reconstruction}, we define
\[
 \Delta L_k^{D,N}
 :=\mathsf R\circ(L_k^{Y/D_0,\log}\mathsf C)-\mathcal A_k^{D,N},
\]
where $L_k^{Y/D_0,\log}$ acts only on the cap coefficients.

\begin{thm}[= Theorem~\ref{thm:degeneration-compatibility}]\label{thm:degeneration-intro}
Let $X\rightsquigarrow X_1\cup_D X_2$ be a complete projective
semistable degeneration in the sense of
Definition~\ref{def:complete-degeneration}, and assume
$\operatorname{Vir}(D)$. The relative Virasoro operators are compatible
with the degeneration formula: for every $k\ge-1$,
\[
 L_k^XZ^X=\Phi_s\bigl(
 \langle L_k^{X_1/D}\mathbf Z^{X_1/D},\mathbf Z^{X_2/D}\rangle_\hbar
 +\langle\mathbf Z^{X_1/D},L_k^{X_2/D}\mathbf Z^{X_2/D}\rangle_\hbar
 \bigr)
\]
where $\Phi_s$ is the specialization defined in
\eqref{eq:general-specialization}.
\end{thm}

\begin{cor}[= Corollary~\ref{cor:degeneration-virasoro}]\label{cor:degeneration-intro}
For a complete projective semistable degeneration
$X\rightsquigarrow X_1\cup_D X_2$, we have
\[
 \operatorname{Vir}(X_1),\quad \operatorname{Vir}(X_2),\quad
 \operatorname{Vir}(D)
 \quad\Longrightarrow\quad \operatorname{Vir}(X).
\]
\end{cor}

\subsection{Virasoro constraints for Enriques surfaces}

For an Enriques surface $S$, the group $H_2(S,\mathbb Z)\cong\mathbb
Z^{10}\oplus\mathbb Z/2$ has torsion generated by the Poincar\'e dual of
the canonical class, and the Gromov--Witten invariants of $S$ are indexed
by integral curve classes.
To recover the invariants of individual integral curve classes from the degeneration
formula, we impose injectivity of specialization on $H_2(S,\mathbb Z)$, namely the
condition (C3) of Definition~\ref{def:complete-degeneration}.
The degeneration of Maulik
and Pandharipande \cite{MP08}*{Section~1.3} kills the torsion class.
Instead, let $R$ be the blow-up of $\PP^2$ at the ten nodes of an
irreducible rational nodal plane sextic, and let $C\subset R$ be
its strict transform. Then $C\cong\PP^1$, $C\sim-2K_R$, and $C^2=-4$.
After contracting $C$, we choose a suitable projective
$\mathbb Q$-Gorenstein smoothing and blow up the singular point
of its total space. As shown in Section~\ref{sec:enriques},
this yields a complete degeneration
\[
 S\rightsquigarrow R\cup_C\PP^2,
\]
in which $C$ is a conic in $\PP^2$.
Corollary~\ref{cor:degeneration-intro} then gives the following result.

\begin{thm}[= Theorem~\ref{thm:enriques}]\label{thm:enriques-intro}
Every Enriques surface $S$ satisfies $\Vir(S)$ in every curve
class $\beta\in H_2(S,\mathbb Z)$.
\end{thm}

\subsection{Organization of the paper}
Section~\ref{sec:relative-virasoro-construction} constructs the relative operators and
proves Theorem~\ref{thm:relative-virasoro-main}.
Section~\ref{sec:examples} contains examples.
Section~\ref{sec:degeneration} proves the compatibility with the
degeneration formula, and Section~\ref{sec:enriques} treats
Enriques surfaces. Appendix~\ref{app:notation-index} collects
principal notations of the main text.

\subsection{Acknowledgment}
We thank Rahul Pandharipande, Pierrick Bousseau, Hülya Argüz and Hiroshi Iritani for helpful discussions. Longting Wu is supported by the National Key R\&D Program of China (No. 2022YFA1006200) and NSFC 12522106.
 Shuai Guo is supported by the NSFC 12225101.

\subsection{Disclosure of AI use}
We used OpenAI’s GPT-6 Astra to 
assist in the construction, develop proof strategies and polish the exposition. We worked out the final mathematical arguments, wrote the manuscript, and take full responsibility for its content and correctness.

\section{Construction of relative Virasoro operators}
\label{sec:relative-virasoro-construction}
In this section, we construct relative Virasoro operators.
We first introduce the gluing
pairing needed to define them. 

\subsection{The gluing pairing}

Write $\eta^{D}_{bc}=\int_D\delta_b\cup\delta_c$
and $(\eta_D^{bc})=(\eta^{D}_{bc})^{-1}$.
We use the left Poincar\'e-dual basis
\[
 \delta^b:=\sum_c\eta_D^{bc}\delta_c,
 \qquad \int_D\delta^b\cup\delta_c
 =\delta_{bc}.
\]
Here $\delta_{bc}$ denotes the Kronecker delta.

We define the gluing pairing on generators $p_{\mu,b}$  by
\begin{equation}\label{eq:gluing-generators}
 \bigl\langle p_{\mu,b},p_{\mu',c}\bigr\rangle_\hbar
 =\delta_{\mu\mu'}\,\mu\hbar\,\eta_D^{bc}.
\end{equation}
Let $p_{\mu,b}^\vee$ denote the coordinate for contact order $\mu$
and the left-dual basis element $\delta^b$. Then
\[
 p_{\mu,b}^\vee
 =\sum_c\eta^D_{cb}p_{\mu,c}.
\]
So
\[\bigl\langle p_{\mu,b},p_{\mu',c}^{\vee}\bigr\rangle_\hbar
 =\delta_{\mu\mu'}\delta_{bc}\,\mu\hbar.\]
For an ordered relative condition
$\lambda=((\mu_j,\delta_{b_j}))_{j=1}^{\ell}$
with no repeated odd generator, define
\[
 \lambda^\vee
 :=((\mu_\ell,\delta^{b_\ell}),\ldots,
    (\mu_1,\delta^{b_1})).
\]
The corresponding coordinate monomial is
\[
 p_\lambda^\vee
 :=p_{\mu_\ell,b_\ell}^\vee\cdots
   p_{\mu_1,b_1}^\vee.
\]
Set $z_\lambda:=|\Aut(\lambda)|\prod_{j=1}^{\ell}\mu_j$.
We extend the gluing pairing bilinearly by
\begin{equation}\label{eq:gluing-pairing}
 \langle p_\lambda,p_{\lambda'}^\vee\rangle_\hbar
 =\delta_{\lambda\lambda'}\,
   \hbar^{\ell(\lambda)}z_\lambda.
\end{equation}
We set $z_{\varnothing}=1$ and $p_{\varnothing}^\vee=1$, so that
$\langle1,1\rangle_\hbar=1$.
All formal variables lie in a common supercommutative coefficient
ring, including variables belonging to different inputs. For
homogeneous coefficients $F,G$ and boundary states $u,v$, we extend
this even pairing by
\[
 \langle Fu,vG\rangle_\hbar
 =F\langle u,v\rangle_\hbar G.
\]
In the gluing pairing, we place coefficients on the left in the first input and on the right in the second. All odd derivatives are left derivatives.

\subsection{Reconstruction via the cap}
\label{subsec:cap-reconstruction}
Let $i:D\hookrightarrow X$ be the inclusion and $N=N_{D/X}$ be the normal bundle of $D$ in $X$. We further set $Y=\PP_D(N\oplus\mathcal O_D)$ to
be the projective bundle over $D$ and let $\pi:Y\to D$ be the projection. Let $D_0$ and $D_\infty$ be the zero and infinity sections of $Y$
such that $N_{D_0/Y}=N^\vee$ and
$N_{D_\infty/Y}=N$. 

Let $[D_0]$ and $[D_\infty]$ be the cohomology classes
of the two sections, and set $\xi=[D_\infty]$ and $\nu=c_1(N)$. It is easy to see that
\[
 [D_0]=\xi-\pi^*\nu,\qquad
 \xi^2=\xi\pi^*\nu,\qquad \xi|_{D_0}=0.
\]

Next, we introduce independent descendant inputs on $D$,
\[
 \mathbf v(\psi)=\sum_{m\ge0,b}v_m^b\delta_b\psi^m,
 \qquad
 \mathbf s(\psi)=\sum_{m\ge0,b}s_m^b\delta_b\psi^m,
\]
and use $\pi^*\mathbf v+\xi\pi^*\mathbf s$ as the ordinary input on $Y$.
The $\mathbf s$-insertions are supported on $D_{\infty}$ and
restrict to zero on $D_0$.

We regard the pair $(Y,D_0)$ as the cap. With the relative conditions placed on the left, define the fixed relative condition
cap potential by
\begin{equation}\label{eq:cap-fixed-potential}
 Z^{Y/D_0}(\lambda\mid\mathbf v,\mathbf s, Q_Y, \hbar)
 :=\sum_{\substack{\beta_Y:\,D_0\cdot\beta_Y=|\lambda|\\g\in\mathbb Z}}
 Q_Y^{\beta_Y}\hbar^{g-1}
 \left\langle\lambda\mathrel{\big|}
 e^{\pi^*\mathbf v+\xi\pi^*\mathbf s}
 \right\rangle^{\bullet,Y/D_0}_{g,\beta_Y}.
\end{equation}
As in \eqref{eq:relative-gw}, the relative markings are labelled,
and the bracket contains no factor of $1/|\Aut(\lambda)|$.
In the defining integral, the relative evaluation classes are
placed before the ordinary descendant insertions.

As in \eqref{eq:relative-total-potential}, we assemble the cap potentials
into the disconnected relative potential:
\begin{equation}\label{eq:cap-state}
\begin{aligned}
 \mathbf Z^{Y/D_0}(\mathbf p;\mathbf v,\mathbf s,Q_Y,\hbar)
 &:=\sum_\lambda
 \frac{p_\lambda}{|\Aut(\lambda)|}
 Z^{Y/D_0}(\lambda\mid\mathbf v,\mathbf s,Q_Y,\hbar)\\
 &=\sum_\lambda
 \frac{p_\lambda^\vee}{|\Aut(\lambda)|}
 Z^{Y/D_0}(\lambda^\vee\mid\mathbf v,\mathbf s,Q_Y,\hbar).
\end{aligned}
\end{equation}
The second equality expresses the same potential in the coordinates
corresponding to the left Poincar\'e-dual basis. The boundary
coordinates are placed to the left of the scalar cap potentials.

Consider now the deformation to the normal cone,
\[
 \mathcal X=\operatorname{Bl}_{D\times\{0\}}(X\times\mathbb A^1)
 \longrightarrow\mathbb A^1,
\]
The general fiber is $X$ and the central fiber is $X\cup_{D=D_0}Y$. The input $\mathbf t$ lifts by pullback from $X$ and restricts
to $\mathbf t$ on the $X$-component and to $\pi^*i^*\mathbf t$
on the $Y$-component.
The strict transform of $D\times\mathbb A^1$ provides a lift of
$i_*\mathbf s$ that restricts to zero on the $X$-component
and to $\xi\pi^*\mathbf s$ on the $Y$-component.
Here $i^*$ and $i_*$ act coefficientwise on descendant inputs.

Jun Li's degeneration formula~\cites{Jun2}, with the above lifts,
takes the form:
\begin{equation}\label{eq:normal-cone-pairing}
 Z^X(\mathbf t+i_*\mathbf s,Q_X,\hbar)
 =\Phi\left\langle
 \mathbf Z^{X/D}(\mathbf t,Q_X,\hbar;\mathbf p),
 \mathbf Z^{Y/D_0}(\mathbf p;i^*\mathbf t,\mathbf s, Q_Y, \hbar)
 \right\rangle_\hbar.
\end{equation}
Here $\Phi$ is defined by
\begin{equation}\label{eq:normal-cone-specialization}
 \Phi\bigl(Q_X^{\beta_X}Q_Y^{\beta_Y}\bigr)
 =Q^{\beta_X+i_*\pi_*\beta_Y}_X.
\end{equation}

It sums over compatible splittings of a fixed total curve class;
the compatibility condition $D\cdot\beta_X=D_0\cdot\beta_Y$
is imposed by the boundary pairing.
Each matched relative condition $\lambda$ contributes the contact and symmetry
factor
\[
 \frac{z_\lambda}{|\Aut(\lambda)|^2}
 =\frac{\prod_j \mu_j}{|\Aut(\lambda)|}
\]
while $\hbar^{\ell(\lambda)}$ supplies the genus shift under gluing.

Next we define the cap operator on the completed boundary state space
by contraction with the cap potential:
\begin{equation}\label{eq:cap-operator}
 \operatorname{Cap}^{\bullet}(\mathbf v, \mathbf s,Q_Y,\hbar)(a)
 :=\left\langle a,
 \mathbf Z^{Y/D_0}(\mathbf p;\mathbf v,\mathbf s, Q_Y,\hbar)
 \right\rangle_\hbar.
\end{equation}
In terms of basis monomials,
\[
 \operatorname{Cap}^{\bullet}(\mathbf v,\mathbf s, Q_Y,\hbar)(p_\lambda)
 =\frac{\hbar^{\ell(\lambda)}z_\lambda}{|\Aut(\lambda)|}
 Z^{Y/D_0}(\lambda^\vee\mid\mathbf v,\mathbf s,Q_Y,\hbar).
\]
Thus \eqref{eq:normal-cone-pairing} can be written as
\begin{equation}\label{eq:normal-cone-cap}
 Z^X(\mathbf t+i_*\mathbf s,Q_X,\hbar)
 =\Phi\left(
 \operatorname{Cap}^{\bullet}(i^*\mathbf t,\mathbf s,Q_Y,\hbar)
 \bigl(\mathbf Z^{X/D}(\mathbf t,Q_X,\hbar;\mathbf p)\bigr)
 \right).
\end{equation}
Here we specialize $\mathbf v=i^*\mathbf t$ in $\operatorname{Cap}^{\bullet}$.

Equation~\eqref{eq:normal-cone-cap} shows that a cap class $\beta_Y$
contributes to the absolute curve class only through its pushforward
$\beta_D:=\pi_*\beta_Y$. This suggests replacing $Q_Y^{\beta_Y}$ by
$Q_D^{\beta_D}$, the Novikov variable of $D$. Let $f$ be the fiber class
and $s_\infty:D\hookrightarrow Y$ the infinity section. For a fixed
relative condition $\lambda$, no curve-class information is lost, since
\begin{equation*}
 \beta_Y=(s_\infty)_*\beta_D+|\lambda|f.
\end{equation*}
After this replacement, cap expressions should be understood as
families indexed by $d=|\lambda|+\int_{\beta_D}c_1(N)$.

For a cohomology-weighted partition $\kappa=((m_i,\delta_{b_i}))_{i=1}^{\ell(\kappa)}$, define
\begin{equation}\label{app:eq:cap-bracket}
\begin{aligned}
 &\left\langle\lambda\mathrel{\big|}-,\kappa
   \right\rangle^{\bullet,D,N}_{g,\beta_D} :=\frac{1}{|\Aut(\kappa)|}
 \bigg\langle\lambda\mathrel{\big|}-,
  \prod_{i=1}^{\ell(\kappa)}
      \xi\pi^*\delta_{b_i}\,\psi_i^{m_i-1}
 \bigg\rangle^{\bullet,Y/D_0}_{g,(s_\infty)_*\beta_D+|\lambda|f},\\[3pt]
 &\left\langle\lambda\mathrel{\big|}-,\kappa
   \right\rangle^{\bullet,D,N}_{\beta_D} :=\sum_g
  \hbar^{g-1+\ell(\lambda)}
 \left\langle\lambda\mathrel{\big|}-,\kappa   \right\rangle^{\bullet,D,N}_{g,\beta_D}.
\end{aligned}
\end{equation}
Here the middle slot contains additional ordinary
insertions ($\varnothing$ for none), and $\kappa$ supplies the ordinary insertions
$\xi\pi^*\delta_{b_i}\psi_i^{m_i-1}$ with symmetry factor $1/|\Aut(\kappa)|$.

We now expand the cap in $\mathbf s$ to express it as a linear map
in monomial bases. For $\kappa$ as above, set
\begin{equation*}
 s_\kappa:=(-1)^{\sum_{i=1}^{\ell}|\delta_{b_i}|}
 s_{m_\ell-1}^{b_\ell}\cdots s_{m_1-1}^{b_1},
 \qquad \ell=\ell(\kappa).
\end{equation*}
Using \eqref{app:eq:cap-bracket}, the cap operator with weights $Q_D^{\beta_D}$ is defined on basis monomials by the following formula
and extended linearly:
\begin{equation}\label{eq:cap-s-expansion}
\begin{aligned}
 &\operatorname{Cap}^{\bullet,D,N}
       (\mathbf v,\mathbf s,Q_D,\hbar)(p_\lambda)\\
 &\quad:=\frac{z_\lambda}{|\Aut(\lambda)|}
   \sum_{\beta_D}Q_D^{\beta_D}
   \sum_\kappa
   \left\langle\lambda^\vee\mathrel{\big|}
       e^{\pi^*\mathbf v},\kappa
   \right\rangle^{\bullet,D,N}_{\beta_D} s_\kappa.
\end{aligned}
\end{equation}
The genus weight is already included
in the bracket. 

Next, we will normalize this map and select its
degree-matched coefficients to obtain an invertible matrix. Its inverse
will give a left inverse of the cap operator, reconstructing $\mathbf Z^{X/D}$
from $Z^X$ and the cap invariants via \eqref{eq:normal-cone-cap}.
Below we omit the arguments $(\mathbf v,\mathbf s,Q_D,\hbar)$
when no confusion can arise.

We set
\begin{equation}\label{eq:cap-empty-factor}
 C_0:=\operatorname{Cap}^{\bullet,D,N}(\mathbf v,\mathbf s,Q_D,\hbar)(1).
\end{equation}
This is the exponential of the connected contribution with empty relative
condition, hence is invertible. The normalized cap
\begin{equation}\label{eq:cap-normalized}
 \operatorname{Cap}^{\bullet,+,D,N}:=C_0^{-1}\operatorname{Cap}^{\bullet,D,N}
\end{equation}
counts maps for which every connected component has a relative marking.
We use the superscript $\bullet,+$ for these invariants and for the
corresponding brackets in \eqref{app:eq:cap-bracket}. In particular,
$\operatorname{Cap}^{\bullet,+,D,N}(1)=1$.
For $\star=\bullet$ or $\bullet,+$, let
$\operatorname{Cap}_d^{\star,D,N}(p_\lambda)$ denote the terms of
$\operatorname{Cap}^{\star,D,N}(p_\lambda)$ with
$|\lambda|+\int_{\beta_D}c_1(N)=d$, equivalently
$D_\infty\cdot\beta_Y=d$.

Let $\mathscr R_D$ be the Novikov ring of $D$ with
coefficients in $\mathbb C((\hbar))[[\mathbf v]]$, completed using an ample
degree on $D$. Coefficient extraction gives the $\mathscr R_D$-linear map
\begin{equation*}
\begin{aligned}
 G:\quad\widehat{\mathcal H}_D\widehat\otimes\mathscr R_D
 &\longrightarrow\widehat{\mathcal H}_D\widehat\otimes\mathscr R_D,\\
 p_\lambda&\longmapsto
   \sum_\kappa 
   [s_\kappa]\operatorname{Cap}_{|\kappa|}^{\bullet,+,D,N}(p_\lambda)p_\kappa.
\end{aligned}
\end{equation*}
Here $[s_\kappa]F=F_\kappa$ for
$F=\sum_\kappa F_\kappa s_\kappa$, with the coefficients
$F_\kappa$ placed on the left.
The tensor product is completed in contact degree, Novikov degree and
the formal variables. A term of class $\beta_D$ can contribute only
when $|\kappa|=|\lambda|+\int_{\beta_D}c_1(N)$.

With the input index first, we write
$G(p_\lambda)=\sum_\kappa G_{\lambda\kappa}p_\kappa$, where
\begin{equation}\label{eq:cap-coefficient-matrix}
\begin{aligned}
 G_{\lambda\kappa}(\mathbf v,Q_D,\hbar)
 &:= [s_\kappa]\operatorname{Cap}_{|\kappa|}^{\bullet,+,D,N}(p_\lambda)\\
 &=\frac{z_\lambda}{|\Aut(\lambda)|}
   \sum_{\substack{\beta_D:\,
           |\kappa|=|\lambda|+\int_{\beta_D}c_1(N)}}
       Q_D^{\beta_D}
   \left\langle\lambda^\vee\mathrel{\big|}
       e^{\pi^*\mathbf v},\kappa
   \right\rangle^{\bullet,+,D,N}_{\beta_D}.
\end{aligned}
\end{equation}
The contact factor is $z_\lambda/|\Aut(\lambda)|=\prod_j\mu_j$.
At $\beta_D=0$, row and column degrees agree, so
\begin{equation*}
 [Q_D^0]G=\operatorname{diag}(G_0,G_1,\ldots),\qquad
 G_d:=\bigl([Q_D^0]G_{\lambda\kappa}\bigr)_{|\lambda|=|\kappa|=d}.
\end{equation*}
Each $G_d$ is a finite square matrix, since $(\mathcal H_D)_d$ is
finite-dimensional.

\begin{lem}\label{lem:cap-leading-invertibility}
$G$ is invertible.
\end{lem}
\begin{proof}
We first prove that $G_d^{(0)}:=G_d|_{\mathbf v=0,\hbar=0}$ is invertible.
For a normalized cap with $c$ connected components, genera $g_j$ and
$\ell$ relative markings, the exponent of $\hbar$ is
$\sum_jg_j-c+\ell\ge0$. Equality holds precisely when every component
has genus zero and one relative marking.
To compute $G_d^{(0)}$, we therefore consider only disconnected invariants
of this type.
By definition, these invariants have $\beta_D=0$, so all curve classes
involved are multiples of the fiber class $f$.

Consider a connected component corresponding to a part
$(\mu,\delta_b)$ of $\lambda$. Its curve class is $\mu f$ and its
relative condition is $(\mu,\delta^b)$ dual to $(\mu,\delta_b)$. Write its ordinary insertions as
$\xi\pi^*\delta_{b_i}\psi_i^{m_i-1}$, with $i$ ranging over the ordinary
markings on this component. We use $\deg_{\mathbb C}$ for half the
cohomological degree. Since
$\deg_{\mathbb C}\delta^b=\dim_{\mathbb C}D-\deg_{\mathbb C}\delta_b$,
the dimension constraint takes the form
\[
 \mu-1+\deg_{\mathbb C}\delta_b
 =\sum_i\bigl(m_i-1+\deg_{\mathbb C}\delta_{b_i}\bigr).
\]
If all ordinary insertions are $\xi$, including the case of no ordinary
markings, the right-hand side vanishes. Since $\mu\ge1$ and
$\deg_{\mathbb C}\delta_b\ge0$, a nonzero contribution requires
$(\mu,\delta_b)=(1,1_D)$, whose dual relative condition is
$(1,[\mathrm{pt}]_D)$. Conversely, for this relative condition the
same equality forces every ordinary insertion to be $\xi$.

For any cohomology-weighted partition $\lambda$, let $b(\lambda)$ be
the number of parts $(1,1_D)$, and set
\[
 \ell_{\mathrm{red}}(\lambda):=\ell(\lambda)-b(\lambda),\qquad
 |\lambda|_{\mathrm{red}}:=|\lambda|-b(\lambda).
\]
These are the number and total contact order of the remaining parts.
For $|\lambda|=|\kappa|=d$, the matrix entry is
\begin{equation}\label{eq:cap-leading-matrix-entry}
 (G_d^{(0)})_{\lambda\kappa}
 =\frac{z_\lambda}{|\Aut(\lambda)|\,|\Aut(\kappa)|}
 \bigg\langle\lambda^\vee\mathrel{\big|}
   \prod_{i=1}^{\ell(\kappa)}
      \xi\pi^*\delta_{b_i}\,\psi_i^{m_i-1}
 \bigg\rangle^{\bullet,+,Y/D_0}_{1-\ell(\lambda),\,df},
\end{equation}
where $\kappa=((m_i,\delta_{b_i}))_{i=1}^{\ell(\kappa)}$.
Each part of $\lambda$ other than $(1,1_D)$ corresponds to a component
requiring at least one insertion other than $\xi$. There are exactly
$\ell_{\mathrm{red}}(\kappa)$ such insertions, so the entry vanishes unless
\[
 \ell_{\mathrm{red}}(\lambda)\le\ell_{\mathrm{red}}(\kappa).
\]
At equality, each of these components carries exactly one insertion
other than $\xi$.

Let $\xi\pi^*\delta_{b_i}\psi^{m_i-1}$ be this insertion on a component
with relative condition $(\mu,\delta^b)$. Since the component maps to a
fiber, all its evaluation maps to $D$ coincide. A nonzero contribution
therefore requires
\[
 \deg_{\mathbb C}\delta^b+\deg_{\mathbb C}\delta_{b_i}
 \le\dim_{\mathbb C}D,
 \qquad\text{or equivalently}\qquad
 \deg_{\mathbb C}\delta_{b_i}\le\deg_{\mathbb C}\delta_b.
\]
The dimension constraint becomes
\[
 \mu-1+\deg_{\mathbb C}\delta_b
 =m_i-1+\deg_{\mathbb C}\delta_{b_i},
\]
so $\mu\le m_i$. Summing over these components yields
$|\lambda|_{\mathrm{red}}\le|\kappa|_{\mathrm{red}}$.
Thus $G_d^{(0)}$ is block triangular when partitions are grouped by
$(\ell_{\mathrm{red}}(\lambda),|\lambda|_{\mathrm{red}})$ and these
groups are ordered lexicographically
(see Example~\ref{ex:p2-cap-matrices}).

At equality in both comparisons, $\mu=m_i$ on each such component.
The fiber calculation in \cite{MP}*{Section~1.4} gives
\begin{equation}\label{eq:cap-one-point-fiber}
 \left\langle(\mu,\delta^b)\mathrel{\big|}
       \xi\pi^*\delta_{b_i}\,\psi^{\mu-1}
 \right\rangle^{\circ,Y/D_0}_{0,\mu f}
 =\frac1{\mu!}\int_D\delta^b\cup\delta_{b_i}
 =\frac{\delta_{b b_i}}{\mu!}.
\end{equation}
Here we use the superscript $\circ$ to indicate that the domain curve is connected.

Equation~\eqref{eq:cap-one-point-fiber} plus the fact that $|\lambda|=|\kappa|$ forces $\lambda=\kappa$ within
each diagonal block. The degree-one components without ordinary markings
contribute $1$, and the divisor equation supplies the nonzero factor
$d^{b(\kappa)}$ for the remaining insertions of $\xi$ when $d>0$.
Thus each diagonal block is diagonal with nonzero diagonal entries.
For $d=0$ the matrix is $(1)$, so $G_d^{(0)}$ is invertible for every $d$.

It follows by formal power-series inversion in $(\mathbf v,\hbar)$ that
each $G_d$ is invertible. Set $A=[Q_D^0]G$ and $R=G-A$. Every term of
$R$ has positive Novikov degree, so $G$ is invertible with
\begin{equation*}
 G^{-1}=\sum_{n\ge0}(-A^{-1}R)^nA^{-1}.
\end{equation*}
\end{proof}

We can now define the cap reconstruction operator. Its input is a
family $\Psi=(\Psi_d)_{d\in\mathbb Z}$ of formal series indexed by an integer degree \(d\).
Regard $C_0$ as the family with entries
$\operatorname{Cap}_d^{\bullet,D,N}(1)$.
Multiplication of families adds their degrees:
$(AB)_d=\sum_{d_1+d_2=d}A_{d_1}B_{d_2}$.
For families whose extracted coefficients define an element of
$\widehat{\mathcal H}_D\widehat\otimes\mathscr R_D$, set
\begin{equation}\label{eq:cap-reconstruction-operator}
 \operatorname{Rec}_{D,N}(\Psi)
 :=\left([s_\kappa](C_0^{-1}\Psi)_{|\kappa|}\right)_\kappa G^{-1}.
\end{equation}
The resulting row vector is identified with an element of
$\widehat{\mathcal H}_D\widehat\otimes\mathscr R_D$.
For any such element $a=\sum_\lambda a_\lambda p_\lambda$,
it is easy to see that
\begin{equation}\label{eq:cap-reconstruction-left-inverse}
 \operatorname{Rec}_{D,N}
 \left((\operatorname{Cap}_d^{\bullet,D,N}(a))_{d\in\mathbb Z}\right)=a.
\end{equation}

\begin{rmk}\label{rmk:cap-absolute-to-relative}
To apply \eqref{eq:cap-reconstruction-left-inverse} to
\eqref{eq:normal-cone-cap}, specialize the coefficients of
$\operatorname{Rec}_{D,N}$ by
$$Q_D^{\beta_D}\mapsto Q_X^{i_*\beta_D},\qquad\mathbf v\mapsto i^*\mathbf t.$$
This preserves the Novikov completion, since an ample class on $X$
restricts to an ample class on $D$. Let $\Psi_d$ be the sum of the terms
of $Z^X(\mathbf t+i_*\mathbf s,Q_X,\hbar)$ with $[D]\cdot\beta=d$.
Then
\[
 \mathbf Z^{X/D}(\mathbf t,Q_X,\hbar;\mathbf p)
 =\operatorname{Rec}_{D,N}\bigl((\Psi_d)_{d\in\mathbb Z}\bigr).
\]
Thus $\operatorname{Rec}_{D,N}$ gives an operator formulation of
Maulik--Pandharipande's absolute-to-relative reconstruction
\cite{MP}*{Section~2}: it recovers the relative potential from the
absolute potential and the cap invariants. The latter can be determined
by the Gromov--Witten theory of $D$ and $c_1(N)$
\cite{MP}*{Theorem~1}.
\end{rmk}

\begin{rmk}\label{rmk:cap-disconnected-divisor}
For disconnected $D$, use the separate contact degrees and the completed
super tensor product of the constructions for its connected components.
\end{rmk}

\subsection{The construction}\label{subsec:relative-construction}
Assume first that $D$ is connected. Let $L_k^{X,\mathbf s}$ denote
\eqref{eq:absolute-virasoro} with its coefficients evaluated at
$\mathbf t+i_*\mathbf s$, all derivatives being taken in $\mathbf t$
with $\mathbf s$ fixed. Thus
\[
 L_k^{X,\mathbf s}Z^X(\mathbf t+i_*\mathbf s)
 =\left.(L_k^XZ^X)(\mathbf u)\right|_{\mathbf u=\mathbf t+i_*\mathbf s},
\]
where $\mathbf u$ is a temporary ordinary descendant input.
Applying this operator to both sides of \eqref{eq:normal-cone-cap} gives
\begin{equation}\label{eq:absolute-action-on-degeneration}
 L_k^{X,\mathbf s}Z^X(\mathbf t+i_*\mathbf s)
 =L_k^{X,\mathbf s}\bigl(\mathsf C(\mathbf Z^{X/D})\bigr).
\end{equation}
For fixed $(D,N)$, we use the following abbreviations throughout the paper:
\[
 \mathsf C:=\operatorname{Cap}^{\bullet,D,N},\qquad
 \mathsf R:=\operatorname{Rec}_{D,N}.
\]
On the right-hand side of \eqref{eq:absolute-action-on-degeneration}, the variables of $ \mathsf C$
are specialized by
\[
 \mathbf v=i^*\mathbf t,\qquad
 Q_D^{\beta_D}\longmapsto Q_X^{i_*\beta_D}.
\]
Thus, the right-hand
side of \eqref{eq:absolute-action-on-degeneration} already includes the
specialization $\Phi$ in \eqref{eq:normal-cone-cap}. 

Our first task is to compute the three-term difference
\begin{equation}\label{eq:relative-defect-start}
 L_k^{X,\mathbf s}\bigl(\mathsf C(\mathbf Z^{X/D})\bigr)
 -\mathsf C\bigl(L_k^{X/D,\log}\mathbf Z^{X/D}\bigr)-(L_k^{Y/D_0,\log}\mathsf C)(\mathbf Z^{X/D}).
\end{equation}
In the last term, $L_k^{Y/D_0,\log}$ acts on $\mathsf C$
before the same specialization of $\mathbf v$ and $Q_D^{\beta_D}$ is applied. We write
\[
 L_k=V_k+B_k+M_k+\delta_{k,0}\varepsilon,
\]
where $V_k$ and $B_k$ are the first- and second-order differential
parts of \eqref{eq:absolute-virasoro}, $M_k$ is its quadratic
multiplication term, and $\varepsilon$ is its scalar term. 
Superscripts will indicate the target and whether the operator is absolute or logarithmic.

\subsubsection{Zero-order terms}
The zero-order contribution to \eqref{eq:relative-defect-start} is
multiplication by
\[
\begin{aligned}
 &M_k^{X,\mathbf s}-M_k^{X/D,\log}-M_k^{Y/D_0,\log}\\
 &\qquad+\delta_{k,0}\bigl(
   \varepsilon^X-\varepsilon^{X/D,\log}
                  -\varepsilon^{Y/D_0,\log}\bigr).
\end{aligned}
\]
We first show that the scalar
\[  \varepsilon^X-\varepsilon^{X/D,\log} -\varepsilon^{Y/D_0,\log}\]
vanishes by lifting $T_X$ to the degeneration family.
Let $f:\mathcal X\to\mathbb A^1$ be the deformation to the normal
cone. Its fibers satisfy $\mathcal X_t\cong X$ for $t\ne0$,
while the central fiber $\mathcal X_0$ is the union of $X$ and $Y$ along $D=D_0$.

The relative
logarithmic tangent bundle
\[
 \mathcal E=\ker\bigl(T_{\mathcal X}(-\log\mathcal X_0)
       \longrightarrow f^*T_{\mathbb A^1}(-\log\{0\})\bigr)
\]
is locally free of rank $d=\dim_{\mathbb C}X$ and
\begin{equation}\label{eq:family-log-tangent}
 \mathcal E|_{\mathcal X_t}=T_{X},\, t\neq 0,\qquad
 \mathcal E|_X=T_X(-\log D),\qquad
 \mathcal E|_Y=T_Y(-\log D_0).
\end{equation}
Since $f$ is proper and flat and its central fiber is the reduced union
$X\cup_D Y$, specialization of the fiber cycle gives
$[\mathcal X_t]\rightsquigarrow[X]+[Y]$. Pairing with the degree-$d$
characteristic class
$(3-d)c_d(\mathcal E)-2c_1(\mathcal E)c_{d-1}(\mathcal E)$
on the total space therefore gives
\[
\begin{aligned}
 &\frac1{48}\int_{\mathcal X_t}
 \bigl((3-d)c_d(\mathcal E|_{\mathcal X_t})
       -2c_1(\mathcal E|_{\mathcal X_t})c_{d-1}(\mathcal E|_{\mathcal X_t})\bigr)\\
 &\quad=\frac1{48}\int_X
 \bigl((3-d)c_d(\mathcal E|_X)
       -2c_1(\mathcal E|_X)c_{d-1}(\mathcal E|_X)\bigr)\\
 &\qquad+\frac1{48}\int_Y
 \bigl((3-d)c_d(\mathcal E|_Y)
       -2c_1(\mathcal E|_Y)c_{d-1}(\mathcal E|_Y)\bigr).
\end{aligned}
\]
By \eqref{eq:absolute-virasoro}, \eqref{eq:logarithmic-scalar},
and \eqref{eq:family-log-tangent}, this is precisely
\begin{equation}\label{eq:logarithmic-scalar-additivity}
 \varepsilon^X=\varepsilon^{X/D,\log}+\varepsilon^{Y/D_0,\log}.
\end{equation}
Hence the
$\delta_{k,0}\varepsilon$ contribution to
\eqref{eq:relative-defect-start} is zero.

Next, we show that the multiplication term
\[M_k^{X,\mathbf s}-M_k^{X/D,\log}-M_k^{Y/D_0,\log}\]
vanishes. We first lift the inputs $\mathbf t$ and $i_*\mathbf s$ to the degeneration family.

Let $p:\mathcal X\to X$ be the collapse.
The input $\mathbf t$ lifts to $p^*\mathbf t$, whose restrictions to
$X$ and $Y$ are $\mathbf t$ and
$\pi^*i^*\mathbf t$, respectively. A lift of $i_*\mathbf s$ is obtained by pulling $\mathbf s$ back
from $D$ to the strict transform of $D\times\mathbb A^1$ and then
pushing forward to $\mathcal X$. This strict transform is disjoint
from $X$ and meets $Y$ along $D_\infty$. Hence the lift restricts
to $0$ on $X$ and to $\xi\pi^*\mathbf s$ on $Y$. Here, we recall that $\xi=[D_\infty]$.

The class $\rho=c_1(X)$ lifts to $c_1(\mathcal E)$, its restrictions to
the two central components are
\begin{equation*}
\begin{aligned}
 \rho_X^{\log}:=c_1(\mathcal E)|_X
 &=c_1(T_X(-\log D))=\rho-[D],\\
 \rho_Y^{\log}:=c_1(\mathcal E)|_Y
 &=c_1(T_Y(-\log D_0))
   =\pi^*i^*\rho+[D_0]
   =\pi^*c_1(D)+\xi.
\end{aligned}
\end{equation*}
We also use $\rho_X^{\log}$ and $\rho_Y^{\log}$ for the corresponding
cup-product operators.

For $k\ge-1$, the quadratic multiplication term $M_k^{X,\mathbf s}$
depends only on the primary inputs $\mathbf t_0$ and $\mathbf s_0$.
By \eqref{eq:absolute-virasoro},
\[
 M_k^{X,\mathbf s}
 =\frac1{2\hbar}\int_X
   (\mathbf t_0+i_*\mathbf s_0)\,
   \rho^{k+1}(\mathbf t_0+i_*\mathbf s_0).
\]
Here the integral denotes the coefficientwise bilinear pairing,
with the formal super coefficients kept in the order of
\eqref{eq:absolute-virasoro}.
Under the lifts above, the primary input restricts to
$\mathbf t_0$ on $X$ and
$\pi^*i^*\mathbf t_0+\xi\pi^*\mathbf s_0$ on $Y$, while
$c_1(\mathcal E)^{k+1}$ restricts to
$(\rho_X^{\log})^{k+1}$ and $(\rho_Y^{\log})^{k+1}$.
Substituting these restrictions into the specialization
$[\mathcal X_t]\rightsquigarrow[X]+[Y]$ gives
\begin{equation*}
\begin{aligned}
 M_k^{X,\mathbf s}
 ={}&\frac1{2\hbar}\int_X
     \mathbf t_0\,(\rho_X^{\log})^{k+1}\mathbf t_0\\
 &+\frac1{2\hbar}\int_Y
   (\pi^*i^*\mathbf t_0+\xi\pi^*\mathbf s_0)\,
   (\rho_Y^{\log})^{k+1}
   (\pi^*i^*\mathbf t_0+\xi\pi^*\mathbf s_0)\\
 ={}&M_k^{X/D,\log}+M_k^{Y/D_0,\log},
 \qquad k\ge-1.
\end{aligned}
\end{equation*}
This includes $k=-1$, when $c_1(\mathcal E)^{k+1}=1$.
Together with \eqref{eq:logarithmic-scalar-additivity}, it proves
that the entire zero-order contribution to
\eqref{eq:relative-defect-start} vanishes.

\subsubsection{First-order terms}
The first-order contribution to \eqref{eq:relative-defect-start} is
\begin{equation}\label{eq:relative-first-order-start}
V_k^{X,\mathbf s}\bigl(\mathsf C(\mathbf Z^{X/D})\bigr)
 -\mathsf C\bigl(V_k^{X/D,\log}\mathbf Z^{X/D}\bigr)
-(V_k^{Y/D_0,\log}\mathsf C)(\mathbf Z^{X/D}).
\end{equation}

Explicitly, \eqref{eq:absolute-virasoro} gives
\[
\begin{aligned}
 V_k^{X,\mathbf s}
 &=\sum_{m\ge0}\sum_{i=0}^{k+1}\sum_{a,b}
   [m+b_a]^k_i(\rho^i)_a{}^b
   \bigl(\widetilde t_m^a+(i_*\mathbf s_m)^a\bigr)
   \partial_{m+k-i,b},\\
 V_k^{X/D,\log}
 &=\sum_{m\ge0}\sum_{i=0}^{k+1}\sum_{a,b}
   [m+b_a]^k_i\bigl((\rho_X^{\log})^i\bigr)_a{}^b
   \widetilde t_m^a\partial_{m+k-i,b}.
\end{aligned}
\]
Here the coefficients $(i_*\mathbf s_m)^a$ are defined by
\[
 i_*\mathbf s(\psi)
 =\sum_{m\ge0,a}(i_*\mathbf s_m)^a\phi_a\psi^m.
\]
The derivatives in the first line act in $\mathbf t$ with
$\mathbf s$ fixed.
The cap operator $V_k^{Y/D_0,\log}$ is given by the same first-order
formula on $H^*(Y)$, using $\rho_Y^{\log}$ and the ordinary input
$\pi^*\mathbf v+\xi\pi^*\mathbf s$, with the dilaton shift by $1_Y$.
It acts only on $\mathsf C$ before the specialization $\mathbf v=i^*\mathbf t$ and
$Q_D^{\beta_D}\mapsto Q_X^{i_*\beta_D}$.

Since our bases are Hodge-homogeneous, we may first fix a value
$b$ and collect the basis elements with $b_a=b$.
Write $\delta_c\in H^{p_c^D,q_c^D}(D)$.
For each $m$, the corresponding part of the shifted input is
\[
\begin{aligned}
 \gamma_{m,b}
 &:=\sum_{a:b_a=b}\widetilde t_m^a\phi_a
   +\sum_{c:\,p_c^D+(3-d)/2=b}
        s_m^c i_*\delta_c\\
 &=\sum_{a:b_a=b}
   \bigl(\widetilde t_m^a+(i_*\mathbf s_m)^a\bigr)\phi_a.
\end{aligned}
\]
The second equality uses the fact that $i_*$ raises Hodge bidegree
by $(1,1)$. For $\eta=\sum_a\eta^a(\mathbf t,\mathbf s)\phi_a$,
we define its action on a formal function $F$ by
\[
 \partial_{\mathbf t_{m'},\eta}F
 :=\sum_a\eta^a(\mathbf t,\mathbf s)
           \frac{\partial F}{\partial t_{m'}^a}.
\]
This notation represents insertion of $\eta\psi^{m'}$
coefficientwise. Then
\[
 V_k^{X,\mathbf s}
 =\sum_{m\ge0}\sum_{i=0}^{k+1}\sum_b
   [m+b]^k_i\,
   \partial_{\mathbf t_{m+k-i},\rho^i\gamma_{m,b}}.
\]
Lift the two sums defining $\gamma_{m,b}$ term by term as above.
Their restrictions to $X$ and $Y$ are, respectively,
\[
\begin{aligned}
 \gamma_{m,b}^X
 &=\sum_{a:b_a=b}\widetilde t_m^a\phi_a,\\
 \gamma_{m,b}^Y
 &=\sum_{a:b_a=b}\widetilde t_m^a\pi^*i^*\phi_a
   +\sum_{c:\,p_c^D+(3-d)/2=b}
        s_m^c\xi\pi^*\delta_c.
\end{aligned}
\]
Multiplication of this lift by $c_1(\mathcal E)^i$ therefore lifts
$\rho^i\gamma_{m,b}$ to a class whose restrictions are
$(\rho_X^{\log})^i\gamma_{m,b}^X$ and
$(\rho_Y^{\log})^i\gamma_{m,b}^Y$.
Applying the degeneration formula with one additional ordinary
marking carrying $(\rho^i\gamma_{m,b})\psi^{m+k-i}$ gives
\[
\begin{aligned}
 \partial_{\mathbf t_{m+k-i},\rho^i\gamma_{m,b}}
      Z^X(\mathbf t+i_*\mathbf s)
 &=\mathsf C\bigl(
   \partial_{\mathbf t_{m+k-i},(\rho_X^{\log})^i\gamma_{m,b}^X}
   \mathbf Z^{X/D}\bigr)\\
 &\quad+\bigl(
   \partial^Y_{m+k-i,(\rho_Y^{\log})^i\gamma_{m,b}^Y}
   \mathsf C\bigr)(\mathbf Z^{X/D}),
\end{aligned}
\]
where $\partial^Y_{m',\eta}$ is the corresponding cap derivative.

Multiplying by $[m+b]^k_i$ and summing over $m,i,b$, the $X$-terms
give $\mathsf C(V_k^{X/D,\log}\mathbf Z^{X/D})$.
The $Y$-terms give
$(V_k^{Y/D_0,\log}\mathsf C)(\mathbf Z^{X/D})$.
Keeping the formal super coefficients in their original order gives
\[
\begin{aligned}
 V_k^{X,\mathbf s}\bigl(\mathsf C(\mathbf Z^{X/D})\bigr)
 &=V_k^{X,\mathbf s}Z^X(\mathbf t+i_*\mathbf s)\\
 &=\mathsf C\bigl(V_k^{X/D,\log}\mathbf Z^{X/D}\bigr)
   +(V_k^{Y/D_0,\log}\mathsf C)(\mathbf Z^{X/D}).
\end{aligned}
\]
This proves that \eqref{eq:relative-first-order-start} vanishes
for every $k\ge-1$.

\subsubsection{Second-order terms}
The second-order contribution to \eqref{eq:relative-defect-start} is
\begin{equation}\label{eq:relative-second-order-start}
 B_k^{X,\mathbf s}\bigl(\mathsf C(\mathbf Z^{X/D})\bigr)
 -\mathsf C\bigl(B_k^{X/D,\log}\mathbf Z^{X/D}\bigr)
 -(B_k^{Y/D_0,\log}\mathsf C)(\mathbf Z^{X/D}).
\end{equation}
Explicitly, \eqref{eq:absolute-virasoro} gives
\[
\begin{aligned}
 B_k^{X,\mathbf s}=B_k^X
 &=\frac{\hbar}{2}\sum_{i=0}^{k-1}\sum_{m=0}^{k-i-1}
   \sum_{a,c}(-1)^{m+1}[-m-b_a]^k_i
   (\rho^i)^{ac}
   \partial_{m,a}\partial_{k-i-1-m,c},\\
 B_k^{X/D,\log}
 &=\frac{\hbar}{2}\sum_{i=0}^{k-1}\sum_{m=0}^{k-i-1}
   \sum_{a,c}(-1)^{m+1}[-m-b_a]^k_i
   \bigl((\rho_X^{\log})^i\bigr)^{ac}
   \partial_{m,a}\partial_{k-i-1-m,c}.
\end{aligned}
\]
Here the indices are raised using the Poincar\'e pairing on $X$,
as in \eqref{eq:absolute-virasoro}. The translation
$\mathbf t\mapsto\mathbf t+i_*\mathbf s$ leaves the constant-coefficient
operator $B_k^X$ unchanged.
The cap operator $B_k^{Y/D_0,\log}$ is given by the same second-order
formula on $H^*(Y)$, using $\rho_Y^{\log}$, the Poincar\'e pairing
on $Y$, and the Hodge grading on $Y$.
It acts only on $\mathsf C$ in the ordinary input
$\pi^*\mathbf v+\xi\pi^*\mathbf s$, before the specialization
$\mathbf v=i^*\mathbf t$ and $Q_D^{\beta_D}\mapsto Q_X^{i_*\beta_D}$.

We fix a value $b$ and sum over $a$ with $b_a=b$.
Let $\phi^a=\sum_c\eta_X^{ac}\phi_c$ be the left Poincar\'e-dual
basis, and put $m'=k-i-1-m$. By the definition of
$(\rho^i)^{ac}$,
\[
 \sum_{a:b_a=b}\sum_c(\rho^i)^{ac}
      \partial_{m,a}\partial_{m',c}
 =\sum_{a:b_a=b}
      \partial_{\mathbf t_m,\phi_a}
      \partial_{\mathbf t_{m'},\rho^i\phi^a}.
\]
Lift $\phi_a$ by $p^*\phi_a$ and $\rho^i\phi^a$ by
$c_1(\mathcal E)^ip^*\phi^a$. Their restrictions to the two
components are
\[
\begin{array}{c|cc}
 & X & Y\\ \hline
 p^*\phi_a & \phi_a & \pi^*i^*\phi_a\\
 c_1(\mathcal E)^ip^*\phi^a
 & (\rho_X^{\log})^i\phi^a
 & (\rho_Y^{\log})^i\pi^*i^*\phi^a.
\end{array}
\]
By definition, $\mathsf C(\mathbf Z^{X/D})$ is the specialized
pairing of $\mathbf Z^{X/D}$ and $\mathbf Z^{Y/D_0}$.
Distributing the two derivatives between its factors by the graded
Leibniz rule, with the componentwise lifts above, gives the four
terms $XX,XY,YX,YY$:
\begin{equation}\label{eq:relative-second-order-four-terms}
\begin{aligned}
 &\sum_{a:b_a=b}\sum_c(\rho^i)^{ac}
   \partial_{m,a}\partial_{m',c}
   \bigl(\mathsf C(\mathbf Z^{X/D})\bigr)\\
 &\quad=\Phi\Biggl\{\sum_{a:b_a=b}\Bigl[
   \bigl\langle
     \partial_{\mathbf t_m,\phi_a}
     \partial_{\mathbf t_{m'},(\rho_X^{\log})^i\phi^a}\mathbf Z^{X/D},
     \mathbf Z^{Y/D_0}
   \bigr\rangle_\hbar\\
 &\qquad\quad+\bigl\langle
     \partial_{\mathbf t_m,\phi_a}\mathbf Z^{X/D},
     \partial^Y_{m',(\rho_Y^{\log})^i\pi^*i^*\phi^a}
       \mathbf Z^{Y/D_0}
   \bigr\rangle_\hbar\\
 &\qquad\quad+(-1)^{|\phi_a|}\bigl\langle
     \partial_{\mathbf t_{m'},(\rho_X^{\log})^i\phi^a}
       \mathbf Z^{X/D},
     \partial^Y_{m,\pi^*i^*\phi_a}\mathbf Z^{Y/D_0}
   \bigr\rangle_\hbar\\
 &\qquad\quad+\bigl\langle
     \mathbf Z^{X/D},
     \partial^Y_{m,\pi^*i^*\phi_a}
     \partial^Y_{m',(\rho_Y^{\log})^i\pi^*i^*\phi^a}\mathbf Z^{Y/D_0}
   \bigr\rangle_\hbar\Bigr]\Biggr\}
   \Bigg|_{\mathbf v=i^*\mathbf t}.
\end{aligned}
\end{equation}
Here $\Phi$ is the Novikov specialization
\eqref{eq:normal-cone-specialization}, $\mathbf v=i^*\mathbf t$
is substituted after differentiation, and $|\phi_a|=p_a+q_a\pmod2$.

We rewrite the last three terms using the diagonal classes of $X$
and $D$. With our left-dual convention,
\[
 [\Delta_X]=\sum_a\phi_a\otimes\phi^a,
 \qquad
 [\Delta_D]=\sum_c\delta_c\otimes\delta^c.
\]
The projection formula makes $i^*$ and $i_*$ adjoint with respect
to the Poincar\'e pairings:
\[
 \int_X u\cup i_*v=\int_D i^*u\cup v,
 \qquad u\in H^*(X),\quad v\in H^*(D).
\]
Expanding this identity in the left-dual bases gives
\[
 (1\otimes i^*)[\Delta_X]=(i_*\otimes1)[\Delta_D],
 \qquad
 (i^*\otimes1)[\Delta_X]=(1\otimes i_*)[\Delta_D].
\]
Recalling that $\nu=c_1(N_{D/X})$, the self-intersection formula
$i^*i_*\gamma=\nu\cup\gamma$ gives
\[
 (i^*\otimes i^*)[\Delta_X]=(\nu\otimes1)[\Delta_D].
\]
Since the bases are Hodge-homogeneous, selecting the summands with
$b_a=b$ gives
\begin{equation*}
\begin{aligned}
 \sum_{a:b_a=b}\phi_a\otimes i^*\phi^a
 &=\sum_{c:\,p_c^D+(3-d)/2=b}
       i_*\delta_c\otimes\delta^c,\\
 \sum_{a:b_a=b}i^*\phi_a\otimes\phi^a
 &=\sum_{c:\,p_c^D+(1-d)/2=b}
       \delta_c\otimes i_*\delta^c,\\
 \sum_{a:b_a=b}i^*\phi_a\otimes i^*\phi^a
 &=\sum_{c:\,p_c^D+(3-d)/2=b}
       \nu\delta_c\otimes\delta^c.
\end{aligned}
\end{equation*}
The shift by one in the first and third sums comes from $i_*$
and multiplication by $\nu$, respectively.

Adjunction and the projection formula also give
\[
(\rho_X^{\log})^ii_*\gamma=i_*\bigl(c_1(D)^i\gamma\bigr).
\]
Consequently, the sum of the last three terms on the right of
\eqref{eq:relative-second-order-four-terms} is
\begin{equation}\label{eq:relative-second-order-D-terms}
\begin{aligned}
 \Phi\Biggl\{&
 \sum_{c:\,p_c^D+(3-d)/2=b}
 \bigl\langle
   \partial_{\mathbf t_m,i_*\delta_c}\mathbf Z^{X/D},
   \partial^Y_{m',(\rho_Y^{\log})^i\pi^*\delta^c}
      \mathbf Z^{Y/D_0}
 \bigr\rangle_\hbar\\
 &+\sum_{c:\,p_c^D+(1-d)/2=b}(-1)^{|\delta_c|}
 \bigl\langle
   \partial_{\mathbf t_{m'},i_*(c_1(D)^i\delta^c)}
      \mathbf Z^{X/D},
   \partial^Y_{m,\pi^*\delta_c}\mathbf Z^{Y/D_0}
 \bigr\rangle_\hbar\\
 &+\sum_{c:\,p_c^D+(3-d)/2=b}
 \bigl\langle
   \mathbf Z^{X/D},
   \partial^Y_{m,\pi^*(\nu\delta_c)}
   \partial^Y_{m',(\rho_Y^{\log})^i\pi^*\delta^c}\mathbf Z^{Y/D_0}
 \bigr\rangle_\hbar\Biggr\}
 \Bigg|_{\mathbf v=i^*\mathbf t}.
\end{aligned}
\end{equation}
Equation~\eqref{eq:relative-second-order-D-terms} shows that we need
to compute terms of the form
$\partial_{\mathbf t_m,i_*\gamma}\mathbf Z^{X/D}$ with
$\gamma\in H^*(D)$. These derivatives can be computed using the
rubber matrix introduced below.

\subsubsection*{The rubber matrix}

Let
\[
 \overline{\mathcal M}^{\bullet,\sim}_{g,n}
 \bigl(Y/(D_0\cup D_\infty),\beta_Y,\lambda,\lambda'\bigr)
\]
denote the moduli space of possibly disconnected rubber stable maps,
with labelled relative conditions $\lambda$ along $D_0$ and $\lambda'$
along $D_\infty$; see \cite{MP}*{Section~1.5.1}.
Writing $\beta_D=\pi_*\beta_Y$, we have
\[
 \beta_Y=(s_\infty)_*\beta_D+|\lambda|f,
 \qquad |\lambda'|-|\lambda|=\nu\cdot\beta_D.
\]
For $\lambda=((\mu_j,\delta_{b_j}))$,
$\lambda'=((\mu'_j,\delta_{c_j}))$, and
$\gamma_i\in H^*(D)$, define
\begin{equation}\label{eq:rubber-invariants}
\begin{aligned}
 &\left\langle
   \lambda\,\middle|\,\prod_{i=1}^n\gamma_i\psi_i^{m_i}
   \,\middle|\,\lambda'
  \right\rangle^{\bullet,\sim,D,N}_{g,\beta_D}\\
 &\quad:=\int_{\left[
   \overline{\mathcal M}^{\bullet,\sim}_{g,n}
   (Y/(D_0\cup D_\infty),\beta_Y,\lambda,\lambda')
   \right]^{\vir}}
   \prod_{j=1}^{\ell(\lambda)}(\ev_j^0)^*\delta_{b_j}\,
   \prod_{i=1}^n\psi_i^{m_i}\ev_i^*\gamma_i\,
   \prod_{j=1}^{\ell(\lambda')}(\ev_j^\infty)^*\delta_{c_j}.
\end{aligned}
\end{equation}
Here all evaluation maps take values in $D$. The factors are ordered as
displayed, and the bracket contains no relative symmetry factors.

For $m\ge0$ and homogeneous
$\gamma\in H^*(D)$, define the rubber matrix by
\begin{equation}\label{eq:rubber-insertion-matrix}
 P_{m,\gamma}^{D,N}(\mathbf v)(p_\lambda)
 :=\frac{(-1)^{|\gamma||p_\lambda|}z_\lambda}{|\Aut(\lambda)|}
 \sum_{\lambda',\beta_D,g}
 \frac{Q_D^{\beta_D}\hbar^{g-1+\ell(\lambda)}}{|\Aut(\lambda')|}\,
 \left\langle\lambda^\vee\,\middle|\,
      e^{\mathbf v} \gamma\psi^m
       \,\middle|\,\lambda'\right\rangle^{\bullet,\sim,D,N}_{g,\beta_D}
 p_{\lambda'}.
\end{equation}
Here $\lambda^\vee$ is the reversed left-dual profile as before, and we
set $P_{m,\gamma}^{D,N}=0$ for $m<0$.
The matrix has parity $|\gamma|$ and is extended graded-linearly with
coefficients placed on the left.

\begin{lem}
\label{lem:rubber-insertion}
For every smooth pair $(X,D)$ with normal bundle $N=N_{D/X}$,
\begin{equation}\label{eq:relative-insertion-equations}
 \partial_{\mathbf t_m,i_*\gamma}\mathbf Z^{X/D}
 =P_{m,\gamma}^{D,N}(i^*\mathbf t)\mathbf Z^{X/D}.
\end{equation}
\end{lem}

\begin{proof}
Apply the degeneration formula \cite{Jun2} to the deformation to
the normal cone of the pair $(X,D)$, using the strict-transform lift
of $i_*\gamma$, which restricts to $0$ on $X$ and to
$\xi\pi^*\gamma$ on $Y$.
The rigidification lemma of Maulik--Pandharipande
\cite{MP}*{Lemma~2, Section~1.5.3} converts the resulting
$(Y,D_0\cup D_\infty)$ invariants with insertion
$\xi\pi^*\gamma\psi^m$ into rubber invariants with insertion
$\gamma\psi^m$, with the ordinary cotangent classes preserved.
Summing over the gluing profiles and the remaining ordinary
insertions gives \eqref{eq:relative-insertion-equations} by
\eqref{eq:rubber-insertion-matrix}, with $\mathbf v=i^*\mathbf t$
and $Q_D^{\beta_D}\mapsto Q_X^{i_*\beta_D}$.
\end{proof}

The graded adjoint of a homogeneous operator $P$ is defined by
\[
 \langle Pu,v\rangle_\hbar
 =(-1)^{|P||u|}\langle u,P^\dagger v\rangle_\hbar,
\]
for homogeneous $u$, and is extended linearly with coefficients
placed on the right.
One can easily check that 
\[
\begin{aligned}
 \bigl(P_{m,\gamma}^{D,N}(\mathbf v)\bigr)^\dagger(p_{\lambda'}^\vee)
 ={}&\frac{z_{\lambda'}}{|\Aut(\lambda')|}
 \sum_{\lambda,\beta_D,g}
 \frac{Q_D^{\beta_D}\hbar^{g-1+\ell(\lambda')}}{|\Aut(\lambda)|}\,
 p_\lambda^\vee
 \left\langle\lambda^\vee\,\middle|\,
 e^{\mathbf v}\gamma\psi^m\,\middle|\,\lambda'
 \right\rangle^{\bullet,\sim,D,N}_{g,\beta_D}.
\end{aligned}
\]

Using \eqref{eq:cap-state} and right linearity, we obtain
\[
\begin{aligned}
 &\bigl(P_{m,\gamma}^{D,N}(\mathbf v)\bigr)^\dagger
       \mathbf Z^{Y/D_0}\\
 &\quad=\sum_{\lambda,\lambda',\beta_D,g}
 \frac{z_{\lambda'} Q_D^{\beta_D}\hbar^{g-1+\ell(\lambda')}}
 {|\Aut(\lambda)|\,|\Aut(\lambda')|^2}\,
 p_\lambda^\vee
 \left\langle\lambda^\vee\,\middle|\,
 e^{\mathbf v}\gamma\psi^m\,\middle|\,\lambda'
 \right\rangle^{\bullet,\sim,D,N}_{g,\beta_D}
 Z^{Y/D_0}((\lambda')^\vee\mid\mathbf v,\mathbf s),
\end{aligned}
\]
with the Novikov specialization induced by collapsing the rubber
component to $D_0$. As in the proof of
Lemma~\ref{lem:rubber-insertion}, degeneration along $D_0\subset Y$
with the strict-transform lift of
$[D_0]\pi^*\gamma$, followed by rigidification, gives
this sum with the rubber component placed first.
Since $[D_0]=\xi-\pi^*\nu$, this proves
\begin{equation}\label{eq:cap-rubber-insertion}
 (\partial_{\mathbf s_m,\gamma}
       -\partial_{\mathbf v_m,\nu\gamma})\mathbf Z^{Y/D_0}
 =\bigl(P_{m,\gamma}^{D,N}(\mathbf v)\bigr)^\dagger
      \mathbf Z^{Y/D_0}.
\end{equation}
We will also write $\partial^Y_{m,[D_0]\pi^*\gamma}$ for
$\partial_{\mathbf s_m,\gamma}-\partial_{\mathbf v_m,\nu\gamma}$.

We now return to \eqref{eq:relative-second-order-D-terms}.
Let $b_c^Y:=p_c^D+(1-d)/2$ be the grading eigenvalue of
$\pi^*\delta_c$ on $Y$. Retain $m'=k-i-1-m$ and abbreviate
$P_{m,\gamma}=P_{m,\gamma}^{D,N}(\mathbf v)$.
By Lemma~\ref{lem:rubber-insertion} and adjunction, that expression is
\begin{equation}\label{eq:relative-second-order-adjoint}
\begin{aligned}
 \Phi\Biggl\{&
 \sum_{c:\,b_c^Y+1=b}
 \Bigl\langle\mathbf Z^{X/D},
 \Bigl(P_{m,\delta_c}^\dagger
       \partial^Y_{m',(\rho_Y^{\log})^i\pi^*\delta^c}
 +\partial_{\mathbf v_m,\nu\delta_c}
       \partial^Y_{m',(\rho_Y^{\log})^i\pi^*\delta^c}\Bigr)
       \mathbf Z^{Y/D_0}\Bigr\rangle_\hbar\\
 &+\sum_{c:\,b_c^Y=b}
 \Bigl\langle\mathbf Z^{X/D},
       (-1)^{|\delta_c|}P_{m',c_1(D)^i\delta^c}^\dagger
       \partial^Y_{m,\pi^*\delta_c}
       \mathbf Z^{Y/D_0}\Bigr\rangle_\hbar
 \Biggr\}\Bigg|_{\mathbf v=i^*\mathbf t}.
\end{aligned}
\end{equation}

With our left-dual convention, the diagonal of $Y$ is
\[
 [\Delta_Y]
 =\sum_c\Bigl(
 [D_\infty]\pi^*\delta_c\otimes\pi^*\delta^c
 +\pi^*\delta_c\otimes[D_0]\pi^*\delta^c\Bigr).
\]
Consequently,
\begin{equation}\label{eq:cap-second-order-diagonal}
\begin{aligned}
 B_k^{Y/D_0,\log}
 =\frac\hbar2\sum_{i=0}^{k-1}\sum_{m=0}^{k-i-1}
 \sum_c(-1)^{m+1}\Bigl(&
 [-m-b_c^Y-1]^k_i
 \partial^Y_{m,\xi\pi^*\delta_c}
 \partial^Y_{m',(\rho_Y^{\log})^i\pi^*\delta^c}\\
 &+[-m-b_c^Y]^k_i
 \partial^Y_{m,\pi^*\delta_c}
 \partial^Y_{m',[D_0]\pi^*(c_1(D)^i\delta^c)}\Bigr).
\end{aligned}
\end{equation}
Here we use the fact that $\rho_Y^{\log}|_{D_0}=c_1(D)$.

In \eqref{eq:relative-second-order-adjoint}, move $P^\dagger$
past the cap derivatives by the graded Leibniz rule.
The cap insertion equation \eqref{eq:cap-rubber-insertion}
and $[D_0]+\pi^*\nu=\xi$
identify the resulting double derivatives, after summation, with
$(B_k^{Y/D_0,\log}\mathsf C)(\mathbf Z^{X/D})$ by
\eqref{eq:cap-second-order-diagonal}.
The remaining contribution is
$-\mathsf C(\mathcal A_k^{D,N}\mathbf Z^{X/D})$,
where the even boundary matrix $\mathcal A_k^{D,N}$ is defined
on basis monomials by
\begin{equation}\label{eq:normal-cone-rubber-correction-ordered}
\begin{aligned}
 \mathcal A_k^{D,N}(p_\lambda)
 &=\frac\hbar2\sum_{i=0}^{k-1}\sum_{m=0}^{k-i-1}\sum_c
 (-1)^{m+1}\\
 &\quad\times\Bigl(
 (-1)^{|\delta_c|}[-m-b_c^Y-1]^k_i
 \partial_{\mathbf v_{m'},c_1(D)^i\delta^c}
 P_{m,\delta_c}\\
 &\qquad\quad+[-m-b_c^Y]^k_i
 \partial_{\mathbf v_m,\delta_c}
 P_{m',c_1(D)^i\delta^c}\Bigr)(p_\lambda).
\end{aligned}
\end{equation}
Here the derivatives act only on the matrix coefficients of $P$.
We extend $\mathcal A_k^{D,N}$ linearly with coefficients on the
left. For $k=-1,0$, we set 
$$\mathcal A_{-1}^{D,N}=\mathcal A_0^{D,N}=0.$$

After summation over $i,m,c$, the two terms in
\eqref{eq:normal-cone-rubber-correction-ordered} give the same contribution.
Expanding $P$ by
\eqref{eq:rubber-insertion-matrix}, the sign from passing the left
derivative through the relative insertion cancels the sign in the
definition of $P$. In the first sum, exchanging the two ordinary
insertions absorbs its factor $(-1)^{|\delta_c|}$, so
\begin{equation}\label{eq:normal-cone-rubber-correction}
\begin{aligned}
 \mathcal A_k^{D,N}(p_\lambda)
 ={}&\frac{z_\lambda}{|\Aut(\lambda)|}
 \sum_{i=0}^{k-1}\sum_{m=0}^{k-i-1}\sum_c
 (-1)^{m+1}[-m-b_c^Y]^k_i\\
 &\quad\cdot\sum_{\lambda',\beta_D,g}
 \frac{Q_D^{\beta_D}\hbar^{g+\ell(\lambda)}}{|\Aut(\lambda')|}\,
 \left\langle\lambda^\vee\,\middle|\,
 e^{\mathbf v}\delta_c\psi_1^m
 \bigl(c_1(D)^i\delta^c\bigr)\psi_2^{m'}
 \,\middle|\,\lambda'\right\rangle^{\bullet,\sim,D,N}_{g,\beta_D}
 p_{\lambda'}.
\end{aligned}
\end{equation}
We have therefore computed \eqref{eq:relative-second-order-start}:
\begin{equation*}
\begin{aligned}
 &B_k^{X,\mathbf s}\bigl(\mathsf C(\mathbf Z^{X/D})\bigr)
 -\mathsf C\bigl(B_k^{X/D,\log}\mathbf Z^{X/D}\bigr)
 -(B_k^{Y/D_0,\log}\mathsf C)(\mathbf Z^{X/D})\\
 &\qquad=-\mathsf C\bigl(\mathcal A_k^{D,N}\mathbf Z^{X/D}\bigr).
\end{aligned}
\end{equation*}
All terms use the specializations in
\eqref{eq:relative-second-order-D-terms}.
Together with the vanishing of the zero- and first-order terms,
this gives
\begin{equation}\label{eq:relative-cap-defect-identity}
 L_k^{X,\mathbf s}\bigl(\mathsf C(\mathbf Z^{X/D})\bigr)
 =\mathsf C\bigl((L_k^{X/D,\log}-\mathcal A_k^{D,N})
                  \mathbf Z^{X/D}\bigr)
  +(L_k^{Y/D_0,\log}\mathsf C)(\mathbf Z^{X/D}).
\end{equation}

We are now ready to construct the relative Virasoro operators.
The first term on the right already contains $L_k^{X/D,\log}$.
We seek a relative operator equal to $L_k^{X/D,\log}$ plus a
boundary correction, so it is natural to apply the reconstruction
operator $\mathsf R$ from \eqref{eq:cap-reconstruction-operator}.
Since $\mathsf R\mathsf C=\mathrm{id}$ by
\eqref{eq:cap-reconstruction-left-inverse}, this removes
$\mathsf C$ from the first term and gives $L_k^{X/D,\log}$
together with the two remaining contributions. We therefore define
\begin{equation}\label{eq:relative-correction-cap-form}
 \Delta L_k^{D,N}
 :=\mathsf R\circ(L_k^{Y/D_0,\log}\mathsf C)-\mathcal A_k^{D,N},
 \qquad
 L_k^{X/D}:=L_k^{X/D,\log}+\Delta L_k^{D,N}(i^*\mathbf t).
\end{equation}
Here $L_k^{Y/D_0,\log}$ differentiates only the cap coefficients,
before the usual specialization of $\mathbf v$ and the Novikov variables.
In the next subsection, we will prove that these operators satisfy
the properties stated in Theorem~\ref{thm:relative-virasoro-main}.

\begin{rmk}\label{rmk:insertion-solution-space}
Let $\mathscr S_{X/D}$ be the space of relative state-valued series
$\mathbf F=\sum_\beta Q_X^\beta\mathbf F_\beta$ such that
$\mathbf F_\beta$ has contact degree $[D]\cdot\beta$ and $\mathbf F$
satisfies
\begin{equation}\label{eq:insertion-solution-space}
 \partial_{\mathbf t_m,i_*\gamma}\mathbf F
 =P_{m,\gamma}^{D,N}(i^*\mathbf t)\mathbf F
 \qquad(m\ge0,\ \gamma\in H^*(D)\text{ homogeneous}).
\end{equation}
In particular,
$\mathbf Z^{X/D}\in\mathscr S_{X/D}$ by
\eqref{eq:relative-insertion-equations}.
A similar calculation of the zero-, first- and second-order terms,
using the graded Leibniz rule and graded adjunction, shows that, for
every $\mathbf F\in\mathscr S_{X/D}$ and $k\ge-1$,
\[
 L_k^{X,\mathbf s}\bigl(\mathsf C(\mathbf F)\bigr)
 =\mathsf C\bigl((L_k^{X/D,\log}-\mathcal A_k^{D,N})\mathbf F\bigr)
  +(L_k^{Y/D_0,\log}\mathsf C)(\mathbf F).
\]
\end{rmk}

\subsection{Proof of Theorem~\ref{thm:relative-virasoro-main}}\label{subsec:relative-proof}

\begin{lem}\label{lem:cap-correction-in-image}
Assume $\operatorname{Vir}(D)$. Then, for every $k\ge-1$,
\begin{equation}\label{eq:cap-correction-in-image}
 L_k^{Y/D_0,\log}\mathsf C-\mathsf C\circ\mathcal A_k^{D,N}= \mathsf C\circ\Delta L_k^{D,N}.
\end{equation}
\end{lem}

\begin{proof}
Since $Y$ is a $\PP^1$-bundle over $D$, $\operatorname{Vir}(Y)$ holds by \cite{CGT}.
For the pair $(Y,D_\infty)$, consider
\[
 \mathbf T:=\mathbf Z^{Y/D_\infty}
 \bigl(\pi^*\mathbf v+[D_0]\pi^*\mathbf w\bigr),
\]
where $\mathbf w$ is a new set of descendant variables, independent of $\mathbf v$.
Applying \eqref{eq:relative-cap-defect-identity} to $(Y,D_\infty)$,
with derivatives taken before restricting the ordinary input, gives
\[
 (L_k^{Y/D_0,\log}\mathsf C)(\mathbf T)
 =\mathsf C\bigl((\mathcal A_k^{D,N}
       -L_k^{Y/D_\infty,\log})\mathbf T\bigr).
\]
Applying $1-\mathsf C\mathsf R$ and using
\eqref{eq:cap-reconstruction-left-inverse}, we obtain
\[
 (1-\mathsf C\mathsf R)
 \bigl((L_k^{Y/D_0,\log}\mathsf C)(\mathbf T)\bigr)=0.
\]
After replacing $Q_Y^{\beta_Y}$ in $\mathbf T$ by
$Q_D^{\pi_*\beta_Y}$, the left-hand side
can be viewed as an application of $\operatorname{Cap}^{\bullet,D,N^\vee}$ in $\mathbf w$,
so applying $\operatorname{Rec}_{D,N^\vee}$
gives, by \eqref{eq:cap-reconstruction-left-inverse},
\[
 (1-\mathsf C\mathsf R)(L_k^{Y/D_0,\log}\mathsf C)=0.
\]
Composing the first identity in \eqref{eq:relative-correction-cap-form}
with $\mathsf C$ and using the preceding identity gives
\eqref{eq:cap-correction-in-image}.
\end{proof}

\begin{proof}[Proof of Theorem~\ref{thm:relative-virasoro-main}]
Assume first that $D$ is connected and let $N=N_{D/X}$.
By \eqref{eq:relative-correction-cap-form},
\[
 L_k^{X/D}=L_k^{X/D,\log}
   +\Delta L_k^{D,N}(i^*\mathbf t).
\]
The correction is local, depending only on $(D,N)$. This proves (1).

For (2) and (3), assume $\operatorname{Vir}(D)$.
Second, combining \eqref{eq:relative-cap-defect-identity},
\eqref{eq:cap-correction-in-image} and \eqref{eq:normal-cone-cap} gives
\begin{equation*}
 L_k^{X,\mathbf s}Z^X(\mathbf t+i_*\mathbf s)
 =\mathsf C\bigl(L_k^{X/D}\mathbf Z^{X/D}\bigr).
\end{equation*}
Using
the injectivity of the cap from
\eqref{eq:cap-reconstruction-left-inverse} gives
\[
 L_k^XZ^X=0
 \quad\Longleftrightarrow\quad
 L_k^{X,\mathbf s}Z^X(\mathbf t+i_*\mathbf s)=0
 \quad\Longleftrightarrow\quad
 L_k^{X/D}\mathbf Z^{X/D}=0.
\]
This proves (3).

Finally, Remark~\ref{rmk:insertion-solution-space} and
Lemma~\ref{lem:cap-correction-in-image} give, for
$\mathbf F\in\mathscr S_{X/D}$,
\begin{equation}\label{eq:relative-cap-intertwining-general}
 L_k^{X,\mathbf s}\bigl(\mathsf C(\mathbf F)\bigr)
 =\mathsf C\bigl(L_k^{X/D}\mathbf F\bigr).
\end{equation}
By \eqref{eq:cap-rubber-insertion} and graded adjunction,
coefficient differentiation gives, before setting $\mathbf v=i^*\mathbf t$,
\[
 \mathsf C\circ P_{m,\gamma}^{D,N}
 =\bigl(\partial_{\mathbf s_m,\gamma}
       -\partial_{\mathbf v_m,\nu\gamma}\bigr)\mathsf C.
\]
Indeed, on a homogeneous boundary state $u$, both sides equal
\[
 (-1)^{|\gamma||u|}\langle u,P_{m,\gamma}^\dagger
 \mathbf Z^{Y/D_0}\rangle_\hbar.
\]
The chain rule therefore gives, for any relative state-valued series
$\mathbf F$ independent of $\mathbf s$,
\begin{equation}\label{eq:cap-insertion-connection}
 \bigl(\partial_{\mathbf t_m,i_*\gamma}
       -\partial_{\mathbf s_m,\gamma}\bigr)\mathsf C(\mathbf F)
 =\mathsf C\bigl((\partial_{\mathbf t_m,i_*\gamma}
       -P_{m,\gamma}^{D,N}(i^*\mathbf t))\mathbf F\bigr).
\end{equation}
Set $\mathbf u=\mathbf t+i_*\mathbf s$. Since $\bigl(\partial_{\mathbf t_m,i_*\gamma}
       -\partial_{\mathbf s_m,\gamma}\bigr)\mathbf u=0$, the operator$\bigl(\partial_{\mathbf t_m,i_*\gamma}
       -\partial_{\mathbf s_m,\gamma}\bigr)$
annihilates
the coefficients of $L_k^{X,\mathbf s}$, which depend on $\mathbf u$.
So it commutes with $L_k^{X,\mathbf s}$.

For $\mathbf F\in\mathscr S_{X/D}$, its insertion equations and
\eqref{eq:cap-insertion-connection} show that
$(\partial_{\mathbf t_m,i_*\gamma}
 -\partial_{\mathbf s_m,\gamma})\mathsf C(\mathbf F)=0$.
Using \eqref{eq:cap-insertion-connection},
\eqref{eq:relative-cap-intertwining-general} and the commutation above,
we obtain
\[
\begin{aligned}
 &\mathsf C\bigl((\partial_{\mathbf t_m,i_*\gamma}
       -P_{m,\gamma}^{D,N}(i^*\mathbf t))L_k^{X/D}\mathbf F\bigr)\\
 &\quad=\bigl(\partial_{\mathbf t_m,i_*\gamma}
       -\partial_{\mathbf s_m,\gamma}\bigr)
       \mathsf C\bigl(L_k^{X/D}\mathbf F\bigr)\\
 &\quad=\bigl(\partial_{\mathbf t_m,i_*\gamma}
       -\partial_{\mathbf s_m,\gamma}\bigr)
       L_k^{X,\mathbf s}\mathsf C(\mathbf F)\\
 &\quad=L_k^{X,\mathbf s}
       \bigl(\partial_{\mathbf t_m,i_*\gamma}
       -\partial_{\mathbf s_m,\gamma}\bigr)\mathsf C(\mathbf F)=0.
\end{aligned}
\]
Injectivity of $\mathsf C$ therefore shows that $L_k^{X/D}\mathbf F$
satisfies the insertion equations.
Thus $L_k^{X/D}$ preserves $\mathscr S_{X/D}$.

Applying \eqref{eq:relative-cap-intertwining-general} twice and using the
absolute commutation relations \eqref{eq:absolute-virasoro-commutator},
we obtain, for $\mathbf F\in\mathscr S_{X/D}$
\[
 \mathsf C\bigl([L_k^{X/D},L_l^{X/D}]\mathbf F\bigr)
 =(k-l)\mathsf C\bigl(L_{k+l}^{X/D}\mathbf F\bigr).
\]
Injectivity of $\mathsf C$ gives (2).

For disconnected $D=\coprod_sD_s$, the normal-cone degeneration has
one cap for each pair $(D_s,N_s)$, and the cap maps act on the
completed graded tensor product with separate contact degrees.
Since the components are disjoint, for $s\ne t$, we have
\[
 i_s^*i_{t*}=0,\qquad
 (i_s^*\otimes i_t^*)[\Delta_X]
 =(i_s^*i_{t*}\otimes1)[\Delta_{D_t}]=0.
\]
The first identity says that an insertion supported on $D_t$ does
not affect the input of the $D_s$ cap; the second makes the mixed
second-order terms between distinct caps vanish. The zero- and
first-order boundary terms also split, since $[D_s][D_t]=0$.
Thus each component contributes its connected-case correction
$\Delta L_k^{D_s,N_s}(i_s^*\mathbf t)$, and these corrections add.
\end{proof}

\subsection{The explicit formulas for $L_{-1}^{X/D}$ and $L_0^{X/D}$}\label{subsec:relative-string-dilaton}

Suppose first that $D$ is connected. For a relative condition
$\lambda=((\mu_j,\delta_{b_j}))_{j=1}^{\ell}$, where
$\delta_{b_j}\in H^{p_{b_j}^D,q_{b_j}^D}(D)$, we set
\[
 b_\lambda=\sum_{j=1}^{\ell}\left(p_{b_j}^D+\frac{1-d}{2}\right),
 \qquad \mathfrak b(p_\lambda)=b_\lambda p_\lambda.
\]
The corrections are
\begin{equation}\label{eq:relative-zero-mode-boundary}
 \Delta L_{-1}^{D,N}=0,\qquad
 \Delta L_0^{D,N}(p_\lambda)=b_\lambda p_\lambda.
\end{equation}

Indeed, $\mathcal A_{-1}^{D,N}=\mathcal A_0^{D,N}=0$, so
\eqref{eq:relative-correction-cap-form} reduces to
$\Delta L_k^{D,N}=\mathsf R\circ(L_k^{Y/D_0,\log}\mathsf C)$
for $k=-1,0$. The relative string, dilaton and divisor equations,
and the Hodge dimension constraint, give
\[
 (L_{-1}^{Y/D_0,\log}\mathsf C)(p_\lambda)=0,\qquad
(L_0^{Y/D_0,\log}\mathsf C)(p_\lambda)
 =b_\lambda\mathsf C(p_\lambda).
\]
For the second identity, the cap has relative profile $\lambda^\vee$.
If $\delta_c\in H^{p_c^D,q_c^D}(D)$, its Poincar\'e dual has type
$(d-1-p_c^D,d-1-q_c^D)$, so its grading on $Y$ is
\[
 (d-1-p_c^D)+\frac{1-d}{2}=-b_c^Y.
\]
Thus the total grading of $\lambda^\vee$ is $-b_\lambda$. The relative
$L_0$ equation gives the negative of this grading as the eigenvalue,
namely $b_\lambda$. Applying $\mathsf R\mathsf C=\mathrm{id}$ proves
the formulas. For disconnected $D=\coprod_sD_s$, a relative condition
is a tuple $\lambda=(\lambda_s)_s$; the corrections add over the
components, giving $\mathfrak b(p_\lambda)=b_\lambda p_\lambda$ with
$b_\lambda=\sum_s b_{\lambda_s}$. In conclusion, we have
\[L_{-1}^{X/D}=L_{-1}^{X/D,\log},\qquad
L_0^{X/D}=L_0^{X/D,\log}+\mathfrak b.
\]

\section{Examples}\label{sec:examples}

\subsection{Curves}\label{subsec:curve-example}
Let $X$ be a curve and $D$ a disjoint union of points. For each point,
$c_1(D)=0$ and $b_c^Y=0$. In
\eqref{eq:normal-cone-rubber-correction}, only $i=0$ contributes, and
for $k\ge1$ we have $[-m]^k_0=\prod_{j=0}^k(j-m)=0$ whenever
$0\le m\le k-1$. For $k=-1,0$ the sum is empty, so
$\mathcal A_k^{\mathrm{pt},N}=0$ for all $k\ge-1$.

The cap is the relative theory of $(\mathbb P^1,\mathrm{pt})$.
By the relative Virasoro constraints for curves
\cite{OP3}*{Theorem~3}, $L_k^{\mathbb P^1/\mathrm{pt},\log}$
annihilates every cap coefficient. Consequently
\[
 \Delta L_k^{\mathrm{pt},N}=0,\qquad
 L_k^{X/D}=L_k^{X/D,\log}\quad(k\ge-1).
\]
These are the relative operators of Okounkov--Pandharipande.

\subsection{The operator \texorpdfstring{$G$ for $(\PP^n,\PP^{n-1})$}{G for (Pn, P(n-1))}}\label{subsec:pn-example}
Let $X=\PP^n$, $D=\PP^{n-1}$, $n\ge2$, and $N=\mathcal O_D(1)$.
We use the notations of Subsection~\ref{subsec:cap-reconstruction},
set $\mathbf v=0$ throughout this subsection, and choose
\[
 h=\nu=c_1(N),\qquad
 \delta_b=h^b,\qquad \delta^b=h^{n-1-b},\qquad
 0\le b\le n-1,\qquad h^n=0.
\]
Let $\lambda=((\mu_j,\delta_{b_j}))_{j=1}^{\ell}$ and
$\kappa=((m_i,\delta_{b'_i}))_{i=1}^r$, with $|\lambda|=|\kappa|=d$.
As in the proof of Lemma~\ref{lem:cap-leading-invertibility}, the leading
matrix $G_d^{(0)}$ can be determined by
\begin{equation}\label{eq:G}
 (G_d^{(0)})_{\lambda\kappa}
 =\sum_{\kappa_1\sqcup\cdots\sqcup\kappa_\ell=\kappa}
   \prod_{j=1}^{\ell}\mu_j
   \left\langle(\mu_j,\delta^{b_j})\mathrel{\big|}\varnothing,\kappa_j
   \right\rangle^{\circ,D,N}_{0,0}.
\end{equation}
Here we decompose \eqref{eq:cap-leading-matrix-entry} into genus-zero
connected invariants (superscript $\circ$) with $\beta_D=0$.
The sum is taken over all ordered decompositions of $\kappa$,
allowing empty $\kappa_j$. The relative markings remain labelled.
The factors $\mu_j$ come from $z_\lambda/|\Aut(\lambda)|=\prod_j\mu_j$, while each bracket includes
$1/|\Aut(\kappa_j)|$.

To compute these connected invariants, define $I_\mu^\alpha$ by the formal expansion
\begin{equation*}
 I_\mu(z)=\frac1{\prod_{j=1}^{\mu}(h+jz)}
 =\sum_{\alpha\ge0}I_\mu^\alpha h^\alpha z^{-\mu-\alpha}.
\end{equation*}
For a single relative condition $(\mu,\delta^b)$ and ordinary partition
$\kappa$, assume that every $m_i\ge2$.
Applying the fiber-class reduction of \cite{MP}*{Section~1.2}
and the fixed-locus analysis of \cite{Wu18}*{Section~2.2}, we find
that the connected invariant
\[
 \left\langle(\mu,\delta^b)\mathrel{\big|}\varnothing,\kappa
 \right\rangle^{\circ,D,N}_{0,0}
\]
vanishes unless
\begin{equation}\label{eq:support}
 \sum_{i=1}^r(m_i+b'_i-1)=\mu+b-1,\qquad
 \sum_{i=1}^r b'_i\le b,\qquad r\le\mu.
\end{equation}
When these conditions hold, summing over the localization contributions in
\cite{Wu18} gives
\begin{equation*}
 \left\langle(\mu,\delta^b)\mathrel{\big|}\varnothing,\kappa
 \right\rangle^{\circ,D,N}_{0,0}
 =\frac{\mu^{r-2}}{|\Aut(\kappa)|}
 [h^{\,b-\sum_i b'_i}](1-h)
 \prod_{i=1}^r
 \left(\sum_{\alpha=0}^{m_i-2}I_{m_i-1-\alpha}^{\alpha}h^\alpha\right).
\end{equation*}
Since $(-1)^\alpha I_\mu^\alpha>0$, all contributions to the indicated coefficient have the same sign, this value is nonzero. Primary insertions ($m_i=1$) can be reduced
by moving their base classes to the relative marking and applying the
divisor equation.

Returning to \eqref{eq:G}, if every $m_i\ge2$, then
\[
 (G_d^{(0)})_{\lambda\kappa}\ne0
 \quad\Longleftrightarrow\quad
 \text{some decomposition satisfies \eqref{eq:support} on every component}.
\]

\begin{example}[{The pair $(\PP^2,\PP^1)$}]\label{ex:p2-cap-matrices}
Use the lexicographic ordering from the proof of
Lemma~\ref{lem:cap-leading-invertibility} for both rows and columns.
For $d=0,1$, the orders are $\varnothing$ and
$((1,\delta_0))\mid((1,\delta_1))$; the full matrices are
\[
 G_0=(1),\qquad
 G_1=\left(\begin{array}{c|c}
 1&0\\\hline
 0&1
 \end{array}\right).
\]
For $d=2$, use
\[
 \begin{array}{c|l}
 (\ell_{\mathrm{red}}(\lambda),|\lambda|_{\mathrm{red}})&\text{Ordered partitions}\\\hline
 (0,0)&((1,\delta_0),(1,\delta_0))\\[2pt]
 (1,1)&((1,\delta_0),(1,\delta_1))\\[2pt]
 (1,2)&((2,\delta_0)),\quad((2,\delta_1))\\[2pt]
 (2,2)&((1,\delta_1),(1,\delta_1))
 \end{array}
\]
and
\[
 G_2=\left(\begin{array}{c|c|cc|c}
 2&0&-\frac{\hbar}{6}&0&0\\\hline
 0&2&-1&-\frac{\hbar}{12}&0\\\hline
 0&0&1&0&0\\
 0&0&0&1&0\\\hline
 0&0&0&0&1
 \end{array}\right),\qquad
 G_2^{(0)}=\left(\begin{array}{c|c|cc|c}
 2&0&0&0&0\\\hline
 0&2&-1&0&0\\\hline
 0&0&1&0&0\\
 0&0&0&1&0\\\hline
 0&0&0&0&1
 \end{array}\right).
\]
\end{example}

\subsection{A Virasoro commutator for \texorpdfstring{$(\PP^2,\PP^1)$}{(P2, P1)}}\label{subsec:p2-virasoro}
Let $X=\PP^2$, $D=\PP^1$, $N=\mathcal O_D(1)$, and $H=[D]$.
Use the bases $\phi_a=H^a$ for $0\le a\le2$ and
$\delta_0=1$, $\delta_1=h=i^*H$, with $\int_Dh=1$.
Write $\Delta L_k=\Delta L_k^{D,N}$ and set
\[
 K=[L_1^{X/D},L_2^{X/D}]+L_3^{X/D},\qquad
 \mathbf T=t_1^2\cdot1,
\]
where $t_1^2$ records $H^2\psi$ and $1=p_\varnothing$ is the
empty boundary state. Let $Q_X$ and $Q_D$ be the Novikov variables corresponding to
the line classes in $X$ and $D$, respectively. We show that $K\ne0$ on the full relative state
space by computing $\left.[Q_X^0]K(\mathbf T)\right|_{\mathbf t=0}$.

Since $\mathbf v=i^*\mathbf t$, we have $v_m^0=t_m^0$ and
$v_m^1=t_m^1$; every $\Delta L_k$ is independent of $t_m^2$.
The logarithmic operators satisfy
$[L_1^{X/D,\log},L_2^{X/D,\log}]=-L_3^{X/D,\log}$, so
\[
 K=[L_1^{X/D,\log},\Delta L_2]+[\Delta L_1,L_2^{X/D,\log}]+[\Delta L_1,\Delta L_2]+\Delta L_3.
\]
The last two terms retain the factor $t_1^2$, as do the
first-order commutators, since $[V,\Delta]\mathbf T=(V\Delta)\mathbf T$ for a
vector field $V$ and a boundary matrix $\Delta$.
The multiplication terms commute with the corrections, and the
second-order part of $L_1^{X/D,\log}$ has only primary derivatives.
Thus only the term with
$\partial_{t_0^0}\partial_{t_1^2}$ in $L_2^{X/D,\log}$ contributes
at $\mathbf t=0$. Its coefficient in \eqref{eq:absolute-virasoro} is
\[
 -\frac\hbar2\left(\frac12\cdot\frac32\cdot\frac52\right)
 +\frac\hbar2\left(-\frac52\right)\left(-\frac32\right)
                  \left(-\frac12\right)=-\frac{15\hbar}{8}.
\]
The Leibniz rule therefore gives
\begin{equation}\label{eq:p2-empty-reduction}
 \left.[Q_X^0]K(\mathbf T)\right|_{\mathbf t=0}
 =\frac{15\hbar}{8}[Q_D^0]
   \left(\left.\partial_{v_0^0}\Delta L_1\right|_{\mathbf v=0}\right)1.
\end{equation}

To compute this unit derivative, put
$Y=\PP_D(N\oplus\mathcal O_D)$ and $\xi=[D_\infty]$,
and also write $h$ for its pullback to $Y$. We have
\[
\begin{gathered}
 c(T_Y(-\log D_0))=(1+2h)(1+\xi),\qquad
 c(T_Y)=(1+2h)(1+2\xi-h),\\
 \xi^2=\xi h,\quad\int_Y\xi h=1.
\end{gathered}
\]
Consequently,
\[
 \varepsilon^{Y/D_0,\log}=\frac{2-2\cdot5}{48}=-\frac16,
 \qquad
 \varepsilon^Y=\frac{4-2\cdot8}{48}=-\frac14.
\]
Hence
\[
\begin{aligned}
 [L_{-1}^{Y/D_0,\log},L_1^{Y/D_0,\log}]
 &=-2L_0^{Y/D_0,\log}
   +2\bigl(\varepsilon^{Y/D_0,\log}-\varepsilon^Y\bigr)\id\\
 &=-2L_0^{Y/D_0,\log}+\frac16\id.
\end{aligned}
\]
Let $C_0=\mathsf C(1)$ be the unnormalized empty cap.
By Subsection~\ref{subsec:relative-string-dilaton},
$(L_{-1}^{Y/D_0,\log}\mathsf C)(p_\lambda)=0$ for every
$\lambda$, and $L_0^{Y/D_0,\log}C_0=0$.
On coefficients independent of $\mathbf s$, the string vector field is
\[
 S_{\mathbf v}=-\partial_{v_0^0}
       +\sum_{m\ge0}\sum_{b=0}^1v_{m+1}^b\partial_{v_m^b}.
\]
Since $\operatorname{Vir}(\PP^1)$ holds, we may apply
\eqref{eq:cap-correction-in-image} to $1$ and then apply
$L_{-1}^{Y/D_0,\log}$. The product rule and the cap string equation give
\[
 \mathsf C\Bigl(S_{\mathbf v}
       \bigl((\Delta L_1+\mathcal A_1^{D,N})1\bigr)\Bigr)
 =[L_{-1}^{Y/D_0,\log},L_1^{Y/D_0,\log}]C_0
 =\frac16 \mathsf C(1).
\]
Applying $\mathsf R\mathsf C=\mathrm{id}$ and setting $\mathbf v=0$ yields
\[
 [Q_D^0]\left(
 \left.\partial_{v_0^0}(\Delta L_1+\mathcal A_1^{D,N})
 \right|_{\mathbf v=0}\right)1=-\frac16\,1.
\]

For the rubber contribution, only $i=m=m'=0$ occurs in
\eqref{eq:normal-cone-rubber-correction}. The gradings
$b_0^Y=-1/2$ and $b_1^Y=1/2$ give coefficients $-3/4$ and $1/4$.
At $Q_D^0$, the empty input profile forces the output profile
and the curve class to be zero. Differentiating in $v_0^0$ gives
the primary insertions $1,1,h$; rigidification
\cite{MP}*{Lemma~2, Section~1.5.3} replaces one unit by $\xi$.
The string equation and dimension constraints leave only the
connected genus-zero, degree-zero three-point invariant, so
\[
 [Q_D^0]\left(
 \left.\partial_{v_0^0}\mathcal A_1^{D,N}\right|_{\mathbf v=0}
 \right)1=\left(-\frac34+\frac14\right)\int_Y\xi h\,1
 =-\frac12\,1.
\]
Subtracting gives
\begin{equation*}
 [Q_D^0]\left(
 \left.\partial_{v_0^0}\Delta L_1\right|_{\mathbf v=0}
 \right)1=\left(-\frac16+\frac12\right)1=\frac13\,1.
\end{equation*}
Combining this with \eqref{eq:p2-empty-reduction}, we obtain
\begin{equation*}
 \left.[Q_X^0]K(\mathbf T)\right|_{\mathbf t=0}
 =\frac{15\hbar}{8}\cdot\frac13\,1=\frac58\hbar\,1\ne0.
\end{equation*}
The test state $\mathbf T$ has curve and contact degree zero, but does not
belong to $\mathscr S_{X/D}$: since $i_*h=H^2$,
\[
 \left.\partial_{\mathbf t_1,i_*h}\mathbf T\right|_{\mathbf t=0}=1,
 \qquad
 \left.P_{1,h}^{D,N}(i^*\mathbf t)\mathbf T\right|_{\mathbf t=0}=0.
\]
For $\mathbf Z^{X/D}\in\mathscr S_{X/D}$, however,
Theorem~\ref{thm:relative-virasoro-main} gives $K\mathbf Z^{X/D}=0$.

\begin{rmk}\label{rmk:relative-low-mode-commutators}
More generally, let $(X,D)$ be a smooth projective pair of $\dim X\geq 2$, with $D\ne\varnothing$. Assume $\operatorname{Vir}(D)$. On the full relative state space,
\[
 K=[L_1^{X/D},L_2^{X/D}]+L_3^{X/D}\ne0.
\]
The relations involving $L_{-1}^{X/D}$ or $L_0^{X/D}$ nevertheless hold:
\begin{equation*}
\begin{aligned}
 [L_{-1}^{X/D},L_k^{X/D}]&=-(k+1)L_{k-1}^{X/D}
       &&(k\ge0),\\
 [L_0^{X/D},L_k^{X/D}]&=-kL_k^{X/D}
       &&(k\ge-1).
\end{aligned}
\end{equation*}
We leave the verification of these claims to the interested reader.
\end{rmk}

\section{Compatibility with degeneration}
\label{sec:degeneration}

In this section, we study the compatibility of relative Virasoro operators with degeneration. 

\subsection{Completeness of the degeneration family}
Let $f:\mathcal X\to B$ be a
projective morphism from a smooth variety to a smooth curve germ
$(B,0)$, smooth over $B\setminus\{0\}$, with central fiber the reduced
normal-crossing divisor $X_1\cup_DX_2$. The components $X_1,X_2$ are
smooth and meet transversally along a smooth divisor $D$. We first assume
that $D$ is connected; see Remark~\ref{rmk:degeneration-disconnected}
for the disconnected case.
Fix a smooth fiber $X$ and write
\[
 N=N_{D/X_1},\quad N_{D/X_2}\cong N^\vee,\quad
 d=\dim_{\mathbb C}X.
\]
The normal bundle relation follows from the normal-crossing family.
Let $\iota_X,\iota_1,\iota_2$ be the inclusions into $\mathcal X$,
and $i_l:D\hookrightarrow X_l$ for $l=1,2$.

We choose a homogeneous basis $\{\Theta_a\}_{a=1}^{M}$ of
$H^*(\mathcal X)$, and let $t_{m,\mathcal X}^a$ be formal variables of parity
$|\Theta_a|$ for $1\leq a\leq M$ and $m\geq0$. For the input
$\mathbf t_{\mathcal X}=\sum_{m,a}t_{m,\mathcal X}^a\Theta_a\psi^m$, set
$\mathbf t_X=\iota_X^*\mathbf t_{\mathcal X}$ and
$\mathbf t_l=\iota_l^*\mathbf t_{\mathcal X}$ for $l=1,2$. Write $Q_{\mathcal X}$
for the Novikov variables of $\mathcal X$ and define
\[
 \operatorname{Sp}_X(Q_X^\beta)=Q_{\mathcal X}^{\iota_{X*}\beta},\qquad
 \operatorname{Sp}_{12}(Q_{X_1}^{\beta_1}Q_{X_2}^{\beta_2})
       =Q_{\mathcal X}^{\iota_{1*}\beta_1+\iota_{2*}\beta_2}.
\]
As in Section~\ref{subsec:relative-construction}, the first relative
state has coefficients on the left and the second has coefficients on
the right. In particular, for the second component we use
\[
 \mathbf Z^{X_2/D}(\mathbf t_2)
 =\sum_\lambda\frac{p_\lambda^\vee}{|\Aut(\lambda)|}
 Z^{X_2/D}(\lambda^\vee\mid\mathbf t_2),
\]
with the relative insertions before the ordinary insertions, as in
\eqref{eq:cap-state}. Jun Li's degeneration formula \cite{Jun2} takes
the form
\begin{equation}\label{eq:general-degeneration}
 \operatorname{Sp}_X Z^X(\mathbf t_X)
 =\operatorname{Sp}_{12}
  \langle\mathbf Z^{X_1/D}(\mathbf t_1),
         \mathbf Z^{X_2/D}(\mathbf t_2)\rangle_\hbar.
\end{equation}
For a matched profile $\lambda$, the pairing supplies
$\prod_j\mu_j/|\Aut(\lambda)|$ and $\hbar^{\ell(\lambda)}$.
The latter records
$g-1=(g_1-1)+(g_2-1)+\ell(\lambda)$.
For the deformation to the normal cone, using the lifts described in
Section~\ref{subsec:relative-construction} and applying the collapse to
the Novikov variables recovers \eqref{eq:normal-cone-pairing}.

Set $H_\oplus=H^*(X_1)\oplus H^*(X_2)$, with the orthogonal-sum
Poincar\'e pairing $\langle-,-\rangle$ and unit $(1,1)$. Define
\begin{equation}\label{eq:matching-cohomology}
\begin{gathered}
 r(\alpha_1,\alpha_2)=i_1^*\alpha_1-i_2^*\alpha_2,\qquad
 J\gamma=(i_{1*}\gamma,-i_{2*}\gamma),\\
 K=\ker r,\qquad \mathcal J=\operatorname{im}J.
\end{gathered}
\end{equation}
Since $i_1^*i_{1*}=c_1(N)\cup$ and
$i_2^*i_{2*}=-c_1(N)\cup$, we have
$\mathcal J\subset K$. The projection formula gives
\begin{equation}\label{eq:matching-projection}
 \langle J\gamma,\alpha\rangle=\int_D\gamma\,(r\alpha).
\end{equation}
Thus $K=\mathcal J^\perp$ and $\mathcal J=K^\perp$. The quotient
\[
 H=K/\mathcal J,\qquad q:K\longrightarrow H,
\]
has a nondegenerate pairing, an induced unit, and a Hodge decomposition.
As in Section~\ref{subsec:relative-construction}, write
$\rho_{X_l}^{\log}=c_1(T_{X_l}(-\log D))\cup$ and set
$\rho_\oplus^{\log}=\rho_{X_1}^{\log}\oplus\rho_{X_2}^{\log}$.
Adjunction gives $i_l^*\rho_{X_l}^{\log}=c_1(D)i_l^*$ and
$\rho_{X_l}^{\log}i_{l*}=i_{l*}c_1(D)$; hence $\rho_\oplus^{\log}$ preserves
$K,\mathcal J$ and induces $\overline\rho^{\log}$ on $H$.
The restriction map
\[
 \Res:H^*(\mathcal X)\longrightarrow H_\oplus,\qquad
 \Res(\Theta_a)=(\iota_1^*\Theta_a,\iota_2^*\Theta_a)
\]
has image $K$, by the retraction onto the central fiber and its
Mayer--Vietoris sequence.

We now recall the Hodge-theoretic description of the two-component degeneration,
following \cite{PS08}.

Let $H^j_{\lim}$ be the limiting mixed Hodge structure, with
$W_\bullet$ denoting the monodromy weight filtration centered
at $j$; see \cite{PS08}*{Lemma--Definition~11.9}.
By \cite{PS08}*{Corollaries~11.23 and~11.41}, Steenbrink's
weight spectral sequence gives, for a two-component degeneration,
\begin{equation}\label{eq:weight-pieces}
\begin{aligned}
 \operatorname{Gr}_j^W H^j_{\lim}&\cong K^j/\mathcal J^j,\\
 \operatorname{Gr}_{j-1}^W H^j_{\lim}
   &\cong\operatorname{coker}
       (r:H_\oplus^{j-1}\longrightarrow H^{j-1}(D)),\\
 \operatorname{Gr}_{j+1}^W H^j_{\lim}
   &\cong\ker
       (J:H^{j-1}(D)(-1)\longrightarrow H_\oplus^{j+1}).
\end{aligned}
\end{equation}

All other weight pieces vanish. The $(-1)$ is the Tate twist of the
Gysin term. The first isomorphism is induced by specialization
$H^j(X_1\cup_DX_2)\to H^j_{\lim}$.
All three are isomorphisms of Hodge structures with the displayed
twists, where the Hodge filtration on the left is the limiting one.

For precision, write $X_0=X_1\cup_DX_2$. The limiting Hodge
filtration $F_\infty^\bullet$ can be obtained from the relative logarithmic
de Rham complex:
\begin{equation*}
 (H^j_{\lim}\otimes\mathbb C,F_\infty^\bullet)
 \cong
 \left(\mathbb H^j\!\left(X_0,
   \Omega^\bullet_{\mathcal X/B}(\log X_0)|_{X_0}\right),F^\bullet\right).
\end{equation*}
Here $F^p$ is induced by forms of degree at least $p$;
see \cite{PS08}*{Corollaries~11.23 and~11.41}. The restriction
maps on the strata have Hodge degree
zero, and the Gysin maps have Hodge degree one, or degree zero
after the Tate twist. Their strictness with respect to the Hodge filtration gives the quotient filtration
\begin{equation*}
\begin{gathered}
 F^pH^j=q(K^j\cap F^pH_\oplus^j),\\
 \mathcal J^j\cap F^pH_\oplus^j
       =J(F^{p-1}H^{j-2}(D)),\\
 (H^j,F^\bullet)\cong
  (\operatorname{Gr}_j^WH^j_{\lim},F_\infty^\bullet).
\end{gathered}
\end{equation*}
Thus passing to $\operatorname{Gr}_F$ commutes with taking the
kernel and quotient defining $H$.

\begin{lem}\label{lem:complete-cohomology}
Assume that $r$ is surjective and $J$ is injective in every degree.
Then monodromy on $H^*(X)$ is trivial, and there is a unique linear
isomorphism $\Psi:H\to H^*(X)$ such that
\[
 (\Psi\circ q\circ\Res)(\Theta_a)=\iota_X^*\Theta_a.
\]
It preserves pairings and units, and
$\Psi\circ\overline\rho^{\log}=c_1(X)\cup\Psi$. Moreover, $\iota_X^*$ is
surjective and
$\iota_{X*}:H_2(X,\mathbb Q)\to H_2(\mathcal X,\mathbb Q)$ is injective.
\end{lem}

\begin{proof}
The assumptions and \eqref{eq:weight-pieces} make the limiting
structure pure. By \cite{PS08}*{Theorem~11.28}, the logarithm of
unipotent monodromy lowers weight by two, so it vanishes.
The first isomorphism in \eqref{eq:weight-pieces} now identifies
$H^j=K^j/\mathcal J^j$ with $H^j_{\lim}$, whose underlying vector
space is identified with $H^j(X)$, defining $\Psi$.
Since this isomorphism is induced by specialization, it gives
$\Psi\circ q\circ\Res=\iota_X^*$.
By the Mayer--Vietoris sequence and the retraction onto the central
fiber, $\operatorname{im}\Res=K$, so $q\circ\Res$ is surjective.
Any other map satisfying this factorization therefore agrees with
$\Psi$ on all of $H$, proving uniqueness.
Since $\Psi$ is an isomorphism, the same factorization also proves
the surjectivity of $\iota_X^*$.
The fibers $X$ and $X_1+X_2$ are homologous; therefore
\[
 \int_X\iota_X^*(\Theta_a\Theta_{a'})
 =\int_{X_1}\iota_1^*(\Theta_a\Theta_{a'})
       +\int_{X_2}\iota_2^*(\Theta_a\Theta_{a'}),
\]
which proves pairing compatibility. The unit extends over the family.
The relative logarithmic tangent bundle restricts to $T_X$ and to
$T_{X_l}(-\log D)$, proving Chern compatibility. Finally, duality
with the surjection in degree two gives the last assertion over
$\mathbb Q$.
\end{proof}

\begin{defn}\label{def:complete-degeneration}
The degeneration is \emph{rationally complete} if it satisfies (C1) and
(C2), and \emph{complete} if it also satisfies (C3):
\begin{enumerate}
\item[(C1)] $r:H_\oplus^j\to H^j(D)$ is surjective and
$J:H^{j-2}(D)\to H_\oplus^j$ is injective for every $j$;
\item[(C2)] the isomorphism $\Psi$ of
Lemma~\ref{lem:complete-cohomology} preserves Hodge types;
\item[(C3)] $\iota_{X*}:H_2(X,\mathbb Z)\to
H_2(\mathcal X,\mathbb Z)$ is injective.
\end{enumerate}
\end{defn}
Condition (C1) is equivalent to trivial monodromy of the degeneration
on $H^*(X,\mathbb Q)$, by \eqref{eq:weight-pieces} and
\cite{PS08}*{Theorem~11.28 and Corollary~11.41}.
Condition (C2) compares the limiting and the chosen smooth-fiber
Hodge decompositions; equivalently, under (C1), it requires
\begin{equation*}
 \Psi(F^pH^j)=F^pH^j(X)\qquad\text{for all }p,j.
\end{equation*}
It is included because our Virasoro operators use the Hodge index
$p$.
(C1) and (C2) together give the Hodge-compatible identification
$\Psi$ on rational cohomology and injectivity of $\iota_{X*}$ on
$H_2(X,\mathbb Q)$ by surjectivity of $\iota_X^*$, hence we call it \emph{rationally complete}.

Assume (C1). Choose a linear lift $\ell$ of $\iota_X^*$
(that is, $\iota_X^*\circ\ell=\mathrm{id}$) with
$\ell(1_X)=1_{\mathcal X}$ and set $s=(s_1,s_2)=\Res\circ\ell$.
Then $\Psi\circ q\circ s=\mathrm{id}$. For
$\mathcal G=\sum Q_{X_1}^{\beta_1}Q_{X_2}^{\beta_2}\mathcal G_{\beta_1,\beta_2}$ define
\begin{equation}\label{eq:general-specialization}
 (\Phi_s \mathcal G)(\mathbf t;Q_X)
 =\sum_\beta Q_X^\beta
  \sum_{\iota_{1*}\beta_1+\iota_{2*}\beta_2=\iota_{X*}\beta}
       \mathcal G_{\beta_1,\beta_2}(s_1\mathbf t,s_2\mathbf t).
\end{equation}
The ample completion controls the compatible effective splittings.
Under (C3), \eqref{eq:general-degeneration} gives
\begin{equation*}
 Z^X=\Phi_s
   \langle\mathbf Z^{X_1/D},\mathbf Z^{X_2/D}\rangle_\hbar.
\end{equation*}

\subsection{Compatibility of Virasoro operators with degeneration}

We now regard $\mathbf t_1,\mathbf t_2$ as independent
descendant inputs on $X_1,X_2$, respectively, and set
\[
 \mathcal G=\langle\mathbf Z^{X_1/D}(\mathbf t_1),
              \mathbf Z^{X_2/D}(\mathbf t_2)\rangle_\hbar,
 \qquad \mathcal F=\operatorname{Sp}_{12}\mathcal G.
\]
We regard $\mathcal F$ as a formal function of $\mathbf t_\oplus=(\mathbf t_1,\mathbf t_2)$.
Restriction to $K$ means that
$\mathbf t_\oplus(\psi)\in K[[\psi]]$, equivalently
$i_1^*\mathbf t_1=i_2^*\mathbf t_2$.

The same degeneration and
rigidification argument as in Lemma~\ref{lem:rubber-insertion} and
\eqref{eq:cap-rubber-insertion} gives
\[
\begin{aligned}
 \partial_{\mathbf t_{1,m},i_{1*}\gamma}\mathbf Z^{X_1/D}
 &=P_{m,\gamma}^{D,N}(i_1^*\mathbf t_1)\mathbf Z^{X_1/D},\\
 \partial_{\mathbf t_{2,m},i_{2*}\gamma}\mathbf Z^{X_2/D}
 &=\bigl(P_{m,\gamma}^{D,N}(i_2^*\mathbf t_2)\bigr)^\dagger
   \mathbf Z^{X_2/D}.
\end{aligned}
\]
Here $\dagger$ is the graded adjoint of
Section~\ref{subsec:relative-construction}. Thus

\begin{equation}\label{eq:matching-direction-derivative}
\begin{aligned}
 \partial_{m,J\gamma}\mathcal F
 &=\operatorname{Sp}_{12}\Bigl(
 \left\langle
 \partial_{\mathbf t_{1,m},i_{1*}\gamma}\mathbf Z^{X_1/D},
 \mathbf Z^{X_2/D}\right\rangle_\hbar-
 \left\langle\mathbf Z^{X_1/D},
 \partial_{\mathbf t_{2,m},i_{2*}\gamma}\mathbf Z^{X_2/D}
 \right\rangle_\hbar\Bigr)\\
 &=\operatorname{Sp}_{12}\left\langle
 \bigl(P_{m,\gamma}^{D,N}(i_1^*\mathbf t_1)
       -P_{m,\gamma}^{D,N}(i_2^*\mathbf t_2)\bigr)
       \mathbf Z^{X_1/D},\mathbf Z^{X_2/D}\right\rangle_\hbar.
\end{aligned}
\end{equation}

\begin{lem}\label{lem:matching-descent}
For every $m,\gamma$,
$(\partial_{m,J\gamma}\mathcal F)|_K=0$.
Consequently $\mathcal F|_K=q^*\overline{\mathcal F}$ for a unique series
$\overline{\mathcal F}$ on $H$. Under the condition of (C1), for every descendant
input $\mathbf t$ on $X$,
\[
 \overline{\mathcal F}(\Psi^{-1}\mathbf t)
 =\operatorname{Sp}_X Z^X(\mathbf t).
\]

\end{lem}

\begin{proof}
On $K$, $i_1^*\mathbf t_1=i_2^*\mathbf t_2$. So
\eqref{eq:matching-direction-derivative} vanishes. Since
$\mathcal J\subset K$, this is exactly
descent to $K/\mathcal J$. Formula~\eqref{eq:general-degeneration}
identifies the descended series with the smooth-fiber potential.
\end{proof}

Now we assume (C1).
The Mayer--Vietoris
sequence identifies $H^*(X_1\cup_DX_2)$ with $K$ by restriction.
The retraction of $\mathcal X$
onto its central fiber also identifies $H^*(X_1\cup_DX_2)$ with
$H^*(\mathcal X)$.
Lemma~\ref{lem:complete-cohomology} identifies $H^*(X)$ with $H=K/\mathcal J$ using $\Psi$.
Under these identifications, choosing a linear map
$\ell:H^*(X)\to H^*(\mathcal X)$ with
$\iota_X^*\circ\ell=\mathrm{id}$ is equivalent to choosing
a complement $\widetilde H$ of $\mathcal J$ in $K$.
The restriction $q|_{\widetilde H}$ identifies $\widetilde H$ with $H$
and preserves the pairing. Hence $\widetilde H$ is nondegenerate,
giving the orthogonal decomposition
\[
 H_\oplus=H^*(X_1)\oplus H^*(X_2)
 =\widetilde H\oplus\widetilde H^\perp.
\]
Choose $\widetilde H$ compatible with the Hodge decomposition.
The following lemma further decomposes $\widetilde H^\perp$.

\begin{lem}\label{lem:matching-witt-decomposition}
Under (C1), for any complement $\widetilde H$ of $\mathcal J$
in $K$ compatible with the Hodge decomposition, there exists
a subspace $T\subset\widetilde H^\perp$ compatible with the
Hodge decomposition such that
\begin{equation*}
\begin{gathered}
 \widetilde H^\perp=\mathcal J\oplus T,\qquad
 \langle T,T\rangle=0,\\
 \langle J\gamma,T_\delta\rangle=\int_D\gamma\delta,
\end{gathered}
\end{equation*}
where $r|_T:T\to H^*(D)$ is an isomorphism and
$T_\delta=(r|_T)^{-1}\delta$.
\end{lem}

\begin{proof}
Since $r$ is surjective and vanishes on $\widetilde H$,
the preceding decomposition shows that $r|_{\widetilde H^\perp}$
is surjective, with kernel $K\cap\widetilde H^\perp=\mathcal J$.
Choose lifts $T_\delta^0\in\widetilde H^\perp$ of
$\delta\in H^*(D)$, linearly in $\delta$ and preserving Hodge types.
Define $u$ by
\[
 \int_Du(\delta)\epsilon=\langle T_\delta^0,T_\epsilon^0\rangle.
\]
Then $u$ has Hodge degree $(-1,-1)$. Set
\[
 T_\delta=T_\delta^0-\tfrac12J(u(\delta)),\qquad
 T=\{T_\delta:\delta\in H^*(D)\}.
\]
Since $rJ=0$, we have $rT_\delta=\delta$. Moreover,
\eqref{eq:matching-projection} and isotropy of $\mathcal J$ give
\[
 \langle T_\delta,T_\epsilon\rangle
 =\bigl(1-\tfrac12-\tfrac12\bigr)
   \langle T_\delta^0,T_\epsilon^0\rangle=0.
\]
Hence $\widetilde H^\perp=\mathcal J\oplus T$, and the remaining
pairing identity follows from \eqref{eq:matching-projection}.
\end{proof}

\begin{prop}\label{prop:diagonal-reduction}
Assume that the degeneration is rationally complete. For every $k\ge-1$,
\begin{equation}\label{eq:diagonal-reduction}
\begin{aligned}
 (\Psi\circ q)^*\operatorname{Sp}_X(L_k^XZ^X)
 ={}&\operatorname{Sp}_{12}\Bigl(
 \langle L_k^{X_1/D,\log}\mathbf Z^{X_1/D},
         \mathbf Z^{X_2/D}\rangle_\hbar\\
 &\quad+\langle\mathbf Z^{X_1/D},
         L_k^{X_2/D,\log}\mathbf Z^{X_2/D}\rangle_\hbar\\
 &\quad-\langle\mathcal A_k^{D,N}(\mathbf v)\mathbf Z^{X_1/D},
         \mathbf Z^{X_2/D}\rangle_\hbar\Bigr)\Big|_K.
\end{aligned}
\end{equation}
Both sides are regarded as formal functions of descendant inputs
$\mathbf t_\oplus(\psi)\in K[[\psi]]$.
Here $(\Psi\circ q)^*$ substitutes $\mathbf t=(\Psi\circ q)(\mathbf t_\oplus)$,
coefficientwise in $\psi$, and
$\mathbf v=i_1^*\mathbf t_1=i_2^*\mathbf t_2$ on $K$.
The correction $\mathcal A_k^{D,N}$ is given by the rubber invariants in
\eqref{eq:normal-cone-rubber-correction} equivalently by
\eqref{eq:normal-cone-rubber-correction-ordered}.
\end{prop}

\begin{proof}
Use the decomposition $L_k=V_k+B_k+M_k+\delta_{k,0}\varepsilon$
from Section~\ref{subsec:relative-construction}, where $V_k,B_k$
are the first- and second-order differential terms, $M_k$ is
quadratic multiplication, and $\delta_{k,0}\varepsilon$ is the
scalar term.
By Lemma~\ref{lem:complete-cohomology} and (C2), $\Psi$
preserves the pairing, unit and Hodge grading, and intertwines
$\overline\rho^{\log}$ with multiplication by $c_1(X)$.
Lemma~\ref{lem:matching-descent} identifies the corresponding
potentials. Thus the $V_k$-terms and $M_k$-terms on both sides of
\eqref{eq:diagonal-reduction} agree.
As in \eqref{eq:logarithmic-scalar-additivity}, specialization of
the relative logarithmic tangent bundle gives
\begin{equation*}
 \varepsilon^X=\varepsilon^{X_1/D,\log}+\varepsilon^{X_2/D,\log}.
\end{equation*}
Hence only the second-order differential terms $B_k$ can contribute
to the difference.

Choose Hodge-homogeneous left-dual bases $\{\phi_a\},\{\phi^a\}$
on $X$ and $\{\delta_c\},\{\delta^c\}$ on $D$.
Fix a decomposition as in Lemma~\ref{lem:matching-witt-decomposition}.
Write $\widetilde\phi_a,\widetilde\phi^a\in\widetilde H$ for the
lifts of $\phi_a,\phi^a$ under $\Psi\circ q$.
By Lemma~\ref{lem:matching-witt-decomposition}, the diagonal
$[\Delta_{H_\oplus}]=[\Delta_{X_1}]+[\Delta_{X_2}]$ decomposes as
\begin{equation*}
 [\Delta_{H_\oplus}]
 =\sum_a\widetilde\phi_a\otimes\widetilde\phi^a
 +\sum_c\bigl(J\delta_c\otimes T_{\delta^c}
                    +T_{\delta_c}\otimes J\delta^c\bigr).
\end{equation*}
The first sum contributes $(\Psi\circ q)^*\operatorname{Sp}_X(B_k^XZ^X)$:
the Hodge gradings agree, the Chern actions agree modulo $\mathcal J$,
and terms with one derivative in $\mathcal J$ and the other tangent
to $K$ vanish by Lemma~\ref{lem:matching-descent}.
Thus only the two mixed sums contribute to the difference.

For $\delta_c\in H^{p_c^D,q_c^D}(D)$, set
$b_c^Y=p_c^D+(1-d)/2$ and $m'=k-i-1-m$;
$J\delta_c$ and $T_{\delta_c}$ have gradings
$b_c^Y+1$ and $b_c^Y$.
Using 
$$(\rho_\oplus^{\log})^iJ\gamma=J(c_1(D)^i\gamma),\qquad
(\rho_\oplus^{\log})^iT_\delta\equiv T_{c_1(D)^i\delta}\pmod K,$$
we obtain
\begin{equation*}
\begin{aligned}
 &\left.\bigl(B_k^{X_1/D,\log}+B_k^{X_2/D,\log}\bigr)
       \mathcal F\right|_K
 -(\Psi\circ q)^*\operatorname{Sp}_X(B_k^XZ^X)\\
 &=\frac\hbar2\sum_{i,m,c}(-1)^{m+1}\Bigl(
 [-m-b_c^Y-1]^k_i
 \partial_{m,J\delta_c}
 \partial_{m',T_{c_1(D)^i\delta^c}}\\
 &\qquad\qquad\qquad+[-m-b_c^Y]^k_i
 \partial_{m,T_{\delta_c}}
 \partial_{m',J(c_1(D)^i\delta^c)}\Bigr)\mathcal F\Big|_K,
\end{aligned}
\end{equation*}
where $0\le i\le k-1$ and $0\le m\le k-i-1$.
To evaluate these mixed derivatives, differentiate
\eqref{eq:matching-direction-derivative} in a $T$-direction, then
\begin{equation*}
 \left.\partial_{m,T_\delta}\partial_{m',J\gamma}\mathcal F\right|_K
 =\operatorname{Sp}_{12}\left\langle
 \bigl(\partial_{\mathbf v_m,\delta}P_{m',\gamma}^{D,N}\bigr)(\mathbf v)
 \mathbf Z^{X_1/D},\mathbf Z^{X_2/D}\right\rangle_\hbar\Big|_K.
\end{equation*}
Here we use the fact that $rT_\delta=\delta$.
Substituting this identity into the two mixed terms gives
\[
 \operatorname{Sp}_{12}\langle
 \mathcal A_k^{D,N}(\mathbf v)\mathbf Z^{X_1/D},
 \mathbf Z^{X_2/D}\rangle_\hbar\Big|_K
\]
by \eqref{eq:normal-cone-rubber-correction-ordered}, equivalently
the rubber-invariant expression \eqref{eq:normal-cone-rubber-correction}.
Moving this correction to the other side proves
\eqref{eq:diagonal-reduction}. For $k=-1,0$, both the second-order
sum and $\mathcal A_k^{D,N}$ vanish.
\end{proof}

By \eqref{eq:relative-correction-cap-form}, writing the right-hand
side of \eqref{eq:diagonal-reduction} in terms of $L_k^{X_1/D}$ and
$L_k^{X_2/D}$, and moving the second boundary correction to the
first factor, leaves the correction
\[
 -\operatorname{Sp}_{12}\Bigl\langle
 \bigl(\Delta L_k^{D,N}+\mathcal A_k^{D,N}
       +(\Delta L_k^{D,N^\vee})^\dagger\bigr)\mathbf Z^{X_1/D},
 \mathbf Z^{X_2/D}\Bigr\rangle_\hbar\Big|_K.
\]
The following lemma shows that this correction vanishes under the condition of
$\operatorname{Vir}(D)$.

\begin{lem}\label{lem:boundary-cancellation}
Assume $\operatorname{Vir}(D)$. For every $k\ge-1$,
\begin{equation*}
 \Delta L_k^{D,N}+\mathcal A_k^{D,N}
       +(\Delta L_k^{D,N^\vee})^\dagger=0,
\end{equation*}
with the adjoint and coefficient conventions as in
Section~\ref{subsec:relative-construction}.
\end{lem}

\begin{proof}
By \cite{CGT} and Theorem~\ref{thm:relative-virasoro-main}, the cap
$(Y,D_0)$ satisfies the relative Virasoro constraints.
Since $N_{D_0/Y}=N^\vee$, pairing these constraints with a fixed
boundary basis vector $u$ and using adjunction gives
\[\begin{aligned}
 0&=\left\langle u,
       \bigl(L_k^{Y/D_0,\log}+\Delta L_k^{D,N^\vee}\bigr)
       \mathbf Z^{Y/D_0}\right\rangle_\hbar\\
  &=(L_k^{Y/D_0,\log}\mathsf C)(u)
       +\mathsf C\bigl((\Delta L_k^{D,N^\vee})^\dagger u\bigr)\\
  &=\mathsf C\Bigl(\bigl(\Delta L_k^{D,N}+\mathcal A_k^{D,N}
       +(\Delta L_k^{D,N^\vee})^\dagger\bigr)u\Bigr).
\end{aligned}\]
Here $L_k^{Y/D_0,\log}$ acts only on the cap coefficients,
and the last equality uses \eqref{eq:cap-correction-in-image}.
Since $u$ is arbitrary and $\mathsf C$ is injective, the assertion follows.
\end{proof}

\begin{thm}\label{thm:rational-degeneration-compatibility}
Assume that the degeneration is rationally complete and that
$\operatorname{Vir}(D)$ holds.
Then, for every $k\ge-1$,
\begin{equation}\label{eq:relative-degeneration-before-specialization}
\begin{aligned}
 (\Psi\circ q)^*\operatorname{Sp}_X(L_k^XZ^X)
 ={}&\operatorname{Sp}_{12}\Bigl(
 \langle L_k^{X_1/D}\mathbf Z^{X_1/D},\mathbf Z^{X_2/D}\rangle_\hbar\\
 &\quad+\langle\mathbf Z^{X_1/D},
 L_k^{X_2/D}\mathbf Z^{X_2/D}\rangle_\hbar\Bigr)\Big|_K.
\end{aligned}
\end{equation}
\end{thm}

\begin{proof}
Combine Proposition~\ref{prop:diagonal-reduction} with
Lemma~\ref{lem:boundary-cancellation}.
\end{proof}

Under the additional condition (C3), this gives the following statement
in the Novikov variables of the smooth fiber.

\begin{thm}\label{thm:degeneration-compatibility}
Assume that the degeneration is complete and that
$\operatorname{Vir}(D)$ holds.
Then, for every $k\ge-1$
\begin{equation}\label{eq:relative-degeneration-compatibility}
 L_k^XZ^X=\Phi_s\bigl(
 \langle L_k^{X_1/D}\mathbf Z^{X_1/D},\mathbf Z^{X_2/D}\rangle_\hbar
 +\langle\mathbf Z^{X_1/D},L_k^{X_2/D}\mathbf Z^{X_2/D}\rangle_\hbar
 \bigr).
\end{equation}
In particular the right-hand side is independent of $s$.
\end{thm}

\begin{proof}
Evaluate \eqref{eq:relative-degeneration-before-specialization}
on $s(\mathbf t)\in K$, using $\Psi\circ q\circ s=\mathrm{id}$,
and extract the coefficients of $Q_{\mathcal X}^{\iota_{X*}\beta}$.
By (C3), this gives
\eqref{eq:relative-degeneration-compatibility}.
\end{proof}

\begin{rmk}\label{rmk:degeneration-disconnected}
For disconnected $D=\coprod_sD_s$, write $i_{l,s}:D_s\hookrightarrow X_l$.
The matching and Gysin maps become
\[
 r(\alpha_1,\alpha_2)
 =\bigl(i_{1,s}^*\alpha_1-i_{2,s}^*\alpha_2\bigr)_s,\qquad
 J((\gamma_s)_s)
 =\sum_s(i_{1,s*}\gamma_s,-i_{2,s*}\gamma_s).
\]
The gluing state space has a separate contact degree for each $D_s$,
and the boundary corrections are summed over the components. Since
$i_{l,s}^*i_{l,t*}=0$ for $s\ne t$, there are no mixed
Gysin--restriction terms between distinct components. The preceding
arguments apply with these maps and sums, with (C1)--(C3) imposed
on the whole degeneration and $\operatorname{Vir}(D_s)$ on each component.
\end{rmk}

\begin{cor}\label{cor:degeneration-virasoro}
For a complete degeneration,
\[
 \operatorname{Vir}(X_1),\quad \operatorname{Vir}(X_2),\quad
 \operatorname{Vir}(D)\quad\Longrightarrow\quad\operatorname{Vir}(X).
\]
\end{cor}

\begin{proof}
Theorem~\ref{thm:relative-virasoro-main} gives
$\operatorname{Vir}(X_1/D)$ and $\operatorname{Vir}(X_2/D)$.
Then \eqref{eq:relative-degeneration-compatibility} gives
$\operatorname{Vir}(X)$.
\end{proof}

\begin{rmk}\label{rmk:rational-completeness}
For a rationally complete degeneration, the preceding identities and
Corollary~\ref{cor:degeneration-virasoro} remain valid with all curve
classes taken modulo torsion and coefficients summed over integral lifts.
Indeed, (C1) makes $\iota_{X*}$ injective modulo torsion, so the same
coefficient extraction applies, where the curve-class equality in
\eqref{eq:general-specialization} is understood in
$H_2(\mathcal X,\mathbb Z)/\mathrm{tors}$.
\end{rmk}

\subsection{A completeness criterion}

\begin{prop}\label{prop:hodge-tate-complete}
If the smooth fiber $X$ is Hodge--Tate, that is,
\[
 H^{p,q}(X)=0\qquad(p\ne q),
\]
then the degeneration is rationally complete.
Furthermore, if $H_2(X,\mathbb Z)$ is also torsion-free, it is complete.
\end{prop}

\begin{proof}
By \cite{PS08}*{Corollary~11.42}, the Hodge--Tate assumption
implies that the monodromy is trivial on rational cohomology, which gives (C1).

Since the limit is pure, \cite{PS08}*{Corollary~11.25} gives
$\dim F^pH^j=\dim F^pH^j(X)$.
The Hodge--Tate assumption therefore implies
\[
 F^pH^j=
 \begin{cases}
 H^j,&2p\le j,\\
 0,&2p>j,
 \end{cases}
\]
and the same formula holds for $H^j(X)$.
Thus $\Psi(F^pH^j)=F^pH^j(X)$, proving (C2).
If $H_2(X,\mathbb Z)$ is torsion-free,
the rational injectivity in Lemma~\ref{lem:complete-cohomology}
implies (C3).
\end{proof}

\section{Virasoro constraints for Enriques surfaces}\label{sec:enriques}

Let $S$ be an Enriques surface. Its cohomology is of Hodge--Tate type,
so every projective semistable degeneration of $S$ is rationally complete
by Proposition~\ref{prop:hodge-tate-complete}.
To obtain Virasoro constraints for integral
classes we still need condition (C3). This matters: $H_2(S,\mathbb
Z)\cong\mathbb Z^{10}\oplus\mathbb Z/2$, with torsion generated by
$\kappa_S=c_1(K_S)\cap[S]$, and the invariants of $\beta$ and
$\beta+\kappa_S$ differ in general \cite{Obe24}*{Appendix~C}. The
degeneration of \cite{MP08}*{Section~1.3} does not satisfy (C3): its two
components have torsion-free $H_2$, and its gluing curve, a fiber of the
rational elliptic surface, is primitive there, so by the Mayer--Vietoris
sequence its central fiber has torsion-free $H_2$, and $\kappa_S$
specializes to zero. We construct a degeneration whose gluing curve is
divisible by two in both components.

Let $B\subset\PP^2$ be an irreducible rational sextic whose only
singularities are ten ordinary nodes, for instance a general rational
plane sextic, let $R$ be the blow-up of $\PP^2$ at the nodes, and let
$C\subset R$ be the strict transform of $B$. Then $R$ is a Coble surface
of K3 type with boundary $C$ \cite{DZ01}: $C\cong\PP^1$,
$C\sim-2K_R$ and $C^2=-4$. Contracting $C$ gives a normal
surface $W$ with one singular point $w_0$, of type $\frac14(1,1)$.

\begin{prop}\label{prop:enriques-degeneration}
There is a projective semistable degeneration $\mathcal X\to\Delta$ over
a disc, with smooth total space, whose smooth fibers are Enriques
surfaces and whose central fiber is
\[
 X_0=R\cup_CP,\qquad P\cong\PP^2,
\]
where $C$ is a conic in $P$.
\end{prop}

\begin{proof}
By \cite{DK25}*{Theorem~9.1.10, Corollary~9.1.11}, $W$ admits a
projective $\mathbb Q$-Gorenstein smoothing to Enriques surfaces. Since
$H^2(W,\Theta_W)=0$ \cite{DK25}*{Section~9.1}, every local deformation
of $w_0$ is induced by a deformation of $W$
\cite{LP07}*{Proposition~2.4}, so the smoothing can be chosen to induce
at $w_0$ the $\mathbb Q$-Gorenstein smoothing $\{xy-z^2=t\}/\{\pm1\}$ of
$\frac14(1,1)=\{xy=z^2\}/\{\pm1\}$ \cite{KSB88}*{Section~3}, where $-1$
acts on $x,y,z$.
Since $W$ has rational singularities and rational resolution $R$,
we have $H^2(W,\mathcal O_W)=H^2(R,\mathcal O_R)=0$.
An ample line bundle on $W$ therefore extends to this globalization
and is relatively ample after shrinking $\Delta$.
We obtain a projective flat family $\mathcal W\to\Delta$
with central fiber $W$ and with Enriques surfaces as smooth fibers,
which near $w_0$ is
\[
 \C^3/\{\pm1\}\longrightarrow\C,\qquad (x,y,z)\longmapsto xy-z^2,
\]
and which is smooth away from $w_0$, since $W$ is.

Let $\mathcal X$ be the blow-up of $\mathcal W$ at $w_0$. The blow-up of
$\C^3/\{\pm1\}$ at the origin is the total space of $\mathcal
O_{\PP^2}(-2)$, where $\PP^2$ has coordinates $[x:y:z]$. On the chart
$x\neq0$, with coordinates $u=x^2$, $v=y/x$, $w=z/x$,
\[
 t=u(v-w^2),
\]
on the chart $z\neq0$, with $u=z^2$, $v=x/z$, $w=y/z$, we have
$t=u(vw-1)$, and the chart $y\neq0$ is symmetric to the first. Hence
$\mathcal X$ is smooth, and its central fiber is reduced, the union of
the exceptional plane $P=\{u=0\}$ and the strict transform of $W$, which
meet transversally along a conic. The strict transform of $W$ is the
blow-up of $W$ at $w_0$, which is the minimal resolution $R$ of $W$, with
exceptional curve $C$.
\end{proof}

\begin{prop}\label{prop:enriques-complete}
The degeneration of Proposition~\ref{prop:enriques-degeneration} is
complete.
\end{prop}

\begin{proof}
Conditions (C1) and (C2) hold by Proposition~\ref{prop:hodge-tate-complete}.
By Lemma~\ref{lem:complete-cohomology}, the kernel of
$\iota_{S*}:H_2(S,\mathbb Z)\to H_2(\mathcal X,\mathbb Z)$ is contained in
the torsion subgroup $\{0,\kappa_S\}$, so (C3) amounts to
$\iota_{S*}\kappa_S\neq0$.

We regard divisor classes on $R$ and $P$ as elements of $H_2$ by
Poincar\'e duality, and let $h$ be the class of a line in $P$, so that
$H_2(P,\mathbb Z)=\mathbb Zh$. The total space retracts onto $X_0=R\cup_CP$, whose
Mayer--Vietoris sequence reads
\[
 H_2(C,\mathbb Z)\longrightarrow H_2(R,\mathbb Z)\oplus\mathbb Zh
 \longrightarrow H_2(X_0,\mathbb Z)\longrightarrow H_1(C,\mathbb Z)=0.
\]
The first map sends $[C]$ to its classes in $R$ and in $P$ with opposite
signs, that is, to $(-2K_R,-2h)$, since $C$ is the anti-bicanonical
curve of $R$ and a conic in $P$. Hence
\[
 H_2(\mathcal X,\mathbb Z)\cong H_2(X_0,\mathbb Z)\cong
 \frac{H_2(R,\mathbb Z)\oplus\mathbb Zh}{\mathbb Z\cdot2(K_R,h)}.
\]
The class $[K_R,h]$ of $(K_R,h)$ satisfies $2[K_R,h]=0$, and it is
nonzero.

Since $S$ is a fiber, $c_1(K_S)$ is the restriction of
$c_1(K_{\mathcal X})$, and $S$ is homologous to $R+P$. By the projection
formula,
\[
 \iota_{S*}\kappa_S=c_1(K_{\mathcal X})\cap[S]
 =c_1(K_{\mathcal X})\cap[R]+c_1(K_{\mathcal X})\cap[P].
\]
Since $R+P=\{t=0\}$ is principal, $\mathcal O_{\mathcal X}(R)|_R\cong
\mathcal O_R(-C)$ and $\mathcal O_{\mathcal X}(P)|_P\cong\mathcal
O_P(-C)$, and adjunction gives $K_{\mathcal X}|_R=K_R+C=-K_R$ and
$K_{\mathcal X}|_P=K_P+C=-h$. Hence $\iota_{S*}\kappa_S=-[K_R,h]
\neq0$.
\end{proof}

\begin{thm}\label{thm:enriques}
Every Enriques surface $S$ satisfies $\Vir(S)$ in every
curve class $\beta\in H_2(S,\mathbb Z)$.
\end{thm}

\begin{proof}
The surface $R$ is the blow-up of $\PP^2$ at ten points, so its big quantum
cohomology is generically semisimple \cite{Bay04}*{Theorem~3.1.1}, as are
those of $P$ and $C$. Hence $R$, $P$ and $C$ satisfy the Virasoro conjecture
\cites{Giv5,T12}. By Proposition~\ref{prop:enriques-complete} and
Corollary~\ref{cor:degeneration-virasoro}, the smooth fibers of the
degeneration of Proposition~\ref{prop:enriques-degeneration} satisfy
the Virasoro constraints in every integral curve class. Since all
Enriques surfaces are Hodge--Tate and deformation equivalent
\cite{BHPV04}*{Chapter~VIII}, the theorem follows by deformation invariance.
\end{proof}

\begin{rmk}
Assuming the genus-two Virasoro constraints, Maulik and Pandharipande
\cite{MP08}*{Proposition~1} derived a recursion for genus-one
invariants in curve classes modulo torsion. Since
localization and Hodge removal used in \cite{MP08}*{Section~3} preserve the integral curve-class
grading,
by Theorem~\ref{thm:enriques}, the same argument now gives an
unconditional recursion in integral curve classes,  It determines the genus-one invariants from those in
integral isotropic classes.
Oberdieck determined the invariants modulo torsion by other
methods \cite{Obe24}.
\end{rmk}

\appendix
\newpage
\section{Index of notation}\label{app:notation-index}


\begingroup\small
\setlength{\tabcolsep}{5pt}
\renewcommand{\arraystretch}{1.06}
\begin{longtable}{@{}>{\raggedright\arraybackslash}p{0.23\textwidth}
                    >{\raggedright\arraybackslash}p{0.59\textwidth}
                    >{\raggedright\arraybackslash}p{0.13\textwidth}@{}}
\toprule
Symbol & Meaning & Reference \\
\midrule
\endfirsthead
\toprule
Symbol & Meaning & Reference \\
\midrule
\endhead
\midrule
\multicolumn{3}{@{}r@{}}{\footnotesize Continued on the next page}\\
\endfoot
\bottomrule
\endlastfoot

\multicolumn{3}{@{}l@{}}{\textbf{Targets and descendant theory}}\\*[2pt]
$X$, $D=\coprod_s D_s$
 & Smooth projective target, and its smooth
   relative divisor with connected components $D_s$.
 & Thm.~\ref{thm:relative-virasoro-main}\\
$i$, $i_s$, $N$, $N_s$
 & Inclusions $D\hookrightarrow X$ and $D_s\hookrightarrow X$ and
   their respective normal line bundles.
 & Sec.~\ref{subsec:cap-reconstruction}, Thm.~\ref{thm:relative-virasoro-main}\\
$\phi_a$, $\delta_b$, $\delta^b$
 & Hodge-homogeneous basis classes on $X$ and $D$, with
   $\phi_a\in H^{p_a,q_a}(X)$, and the left
   Poincar\'e-dual class to $\delta_b$.
 & \eqref{eq:absolute-virasoro}, \eqref{eq:gluing-generators}\\
$\eta^X_{ab}$, $\eta_X^{ab}$\newline
$\eta^D_{bc}$, $\eta_D^{bc}$
 & Poincar\'e pairing matrices and their inverses on $X$ and $D$.
 & \eqref{eq:absolute-virasoro}, \eqref{eq:gluing-generators}\\
$\overline{\mathcal M}_{g,n}(X,\beta)$\newline
$\overline{\mathcal M}^{\bullet}_{g,n}(X/D,\beta,\mu)$
 & Moduli spaces of absolute and possibly disconnected relative
   stable maps. In the relative space $\mu$ records the contact orders
   of the labelled relative markings.
 & \eqref{eq:absolute-gw}, \eqref{eq:relative-gw}\\
$\psi_i$, $\ev_i$, $\ev_j^{\mathrm{rel}}$\newline
$\gamma\psi^m$
 & Cotangent line class at an ordinary marking, ordinary and relative
   evaluation maps, and the descendant insertion $\ev_i^*\gamma\,\psi_i^m$.
 & \eqref{eq:absolute-gw}, \eqref{eq:rubber-invariants}\\
$\mathbf t(\psi)$, $t_m^a$
 & Ordinary descendant input and its coordinate for $\phi_a\psi^m$.
 & \eqref{eq:absolute-connected-potential}\\
$F_g^X$, $F^X$, $Z^X$
 & Genus-$g$ and all-genus connected absolute descendant potentials
   and the disconnected potential $Z^X=\exp(F^X)$.
 & \eqref{eq:absolute-connected-potential}, \eqref{eq:absolute-total-potential}\\[4pt]

\multicolumn{3}{@{}l@{}}{\textbf{Relative conditions and gluing}}\\*[2pt]
$\langle-\rangle^X_{g,\beta}$\newline
$\langle-\mid\lambda\rangle^{\bullet,X/D}_{g,\beta}$
 & Absolute and relative invariant brackets. The relative bracket
   uses labelled markings with no factor of $1/|\Aut(\lambda)|$.
 & \eqref{eq:absolute-gw}, \eqref{eq:relative-gw}\\
$\lambda$, $\mu_j$
 & Relative profile, a partition with cohomology weights
   $(\mu_j,\delta_{b_j})$. Each $\mu_j$ is a positive contact order.
 & \eqref{eq:relative-gw}\\
$|\lambda|$, $\ell(\lambda)$
 & Total contact order $\sum_j\mu_j$ and number of relative markings.
 & \eqref{eq:relative-gw}\\
$|\Aut(\lambda)|$, $z_\lambda$
 & Symmetry order of the weighted profile and the gluing weight
   $z_\lambda=|\Aut(\lambda)|\prod_j\mu_j$.
 & \eqref{eq:gluing-pairing}\\
$p_{\mu,b}$, $p_\lambda$, $\mathbf p$
 & Boundary generator of contact order $\mu$ and weight $\delta_b$,
   ordered profile monomial, and collection of boundary variables.
 & \eqref{eq:relative-total-potential}\\
$\mathcal H_D$, $\widehat{\mathcal H}_D$
 & Supercommutative boundary state space and its completion by
   contact degree.
 & \eqref{eq:relative-total-potential}\\
$\lambda^\vee$, $p_\lambda^\vee$
 & Reversed profile with left-dual cohomology weights and its
   corresponding dual coordinate monomial.
 & \eqref{eq:gluing-pairing}\\
$Z^{X/D}(\cdot\mid\lambda)$, $\mathbf Z^{X/D}$
 & Disconnected relative potentials.
 & \eqref{eq:relative-disconnected-potential}, \eqref{eq:relative-total-potential}\\
$\langle-,-\rangle_\hbar$
 & Gluing pairing of boundary states.
 & \eqref{eq:gluing-pairing}\\
$\circ$, $\bullet$, $\bullet,+$, $\bullet,\sim$
 & Connected, possibly disconnected, possibly disconnected with a
   relative marking on each connected component, and possibly
   disconnected rubber theory, respectively.
 & \eqref{eq:relative-gw}, \eqref{eq:cap-normalized}, \eqref{eq:rubber-invariants}\\[4pt]

\multicolumn{3}{@{}l@{}}{\textbf{Cap and rubber construction}}\\*[2pt]
$Y$, $\pi$, $D_0$, $D_\infty$
 & Ruled target $\PP_D(N\oplus\mathcal O_D)$, its projection to $D$,
   and its sections with normal bundles $N^\vee$ and $N$.
 & Sec.~\ref{subsec:cap-reconstruction}\\
$\xi$, $\nu$, $\beta_D$
 & Class $[D_\infty]$, first Chern class $c_1(N)$, and
   base curve class $\pi_*\beta_Y$.
 & Sec.~\ref{subsec:cap-reconstruction}\\
$\mathbf v$, $\mathbf s$
 & Descendant inputs on $D$ giving $\pi^*\mathbf v+\xi\pi^*\mathbf s$
   on $Y$. The specialization in the cap is $\mathbf v=i^*\mathbf t$.
 & \eqref{eq:cap-fixed-potential}, \eqref{eq:normal-cone-cap}\\
$\mathbf Z^{Y/D_0}$, $\mathsf C$
 & Boundary state of the cap $(Y,D_0)$ and contraction with that
   state. In the construction $\mathsf C=\operatorname{Cap}^{\bullet,D,N}$
   uses base Novikov variables.
 & \eqref{eq:cap-state}, \eqref{eq:cap-operator}, \eqref{eq:cap-s-expansion}\\
$C_0$, $\operatorname{Cap}^{\bullet,+,D,N}$
 & Empty-profile cap factor and normalized cap $C_0^{-1}\mathsf C$.
 & \eqref{eq:cap-empty-factor}, \eqref{eq:cap-normalized}\\
$G$, $G_{\lambda\kappa}$
 & Matrix of selected normalized cap coefficients, with input
   profile $\lambda$ and output monomial indexed by $\kappa$.
 & \eqref{eq:cap-coefficient-matrix}\\
$G_d$, $G_d^{(0)}$
 & Vertical degree-$d$ block and its leading term:\newline
   $G_d=\bigl([Q_D^0]G_{\lambda\kappa}\bigr)_{|\lambda|=|\kappa|=d}$,\newline
   $G_d^{(0)}=G_d\big|_{\mathbf v=0,\,\hbar=0}$.
 & Sec.~\ref{subsec:cap-reconstruction}\\
$\mathsf R=\operatorname{Rec}_{D,N}$
 & Reconstruction operator and left inverse of $\mathsf C$ on the
   families of cap coefficients graded by $D_\infty\cdot\beta_Y$.
 & \eqref{eq:cap-reconstruction-operator}, \eqref{eq:cap-reconstruction-left-inverse}\\
$P_{m,\gamma}^{D,N}$, $P^\dagger$
 & Rubber insertion matrix for $\gamma\psi^m$ and the graded
   adjoint with respect to $\langle-,-\rangle_\hbar$.
 & \eqref{eq:rubber-insertion-matrix}, \eqref{eq:cap-rubber-insertion}\\[4pt]

\multicolumn{3}{@{}l@{}}{\textbf{Virasoro operators}}\\*[2pt]
$b_a$, $\rho$, $[x]^k_i$
 & The Hodge weight $p_a+(1-d)/2$, cup product with $c_1(X)$, and
   the elementary symmetric coefficients of the Virasoro formula.
 & \eqref{eq:virasoro-coefficients}, \eqref{eq:absolute-virasoro}\\
$\widetilde t_m^a$, $\partial_{m,a}$
 & Dilaton-shifted coordinate and left derivative in $t_m^a$.
 & \eqref{eq:absolute-virasoro}\\
$T_X^{\log}$, $L_k^{X/D,\log}$
 & Logarithmic tangent bundle $T_X(-\log D)$ and the operator
   obtained by using its Chern classes in $L_k^X$.
 & \eqref{eq:logarithmic-scalar}, Thm.~\ref{thm:relative-virasoro-main}\\
$L_k^X$, $L_k^{X/D}$
 & Absolute and relative Virasoro operators.
 & \eqref{eq:absolute-virasoro}, Thm.~\ref{thm:relative-virasoro-main}\\
$\mathcal A_k^{D,N}$, $\Delta L_k^{D,N}$
 & Rubber contribution to the comparison with the cap and the
   correction added to the logarithmic operator.
 & \eqref{eq:normal-cone-rubber-correction}, \eqref{eq:relative-correction-cap-form}\\
$\mathfrak b$, $b_\lambda$
 & Boundary grading operator in $L_0^{X/D}$ and its eigenvalue
   on $p_\lambda$.
 & \eqref{eq:relative-zero-mode-boundary}\\
$\mathscr S_{X/D}$
 & Space of boundary-valued series satisfying the contact-degree
   condition and the rubber insertion equations.
 & \eqref{eq:insertion-solution-space}\\
$\Vir(X)$, $\Vir(X/D)$
 & Assertions that the absolute or relative descendant potential,
   respectively, is annihilated by every corresponding Virasoro operators.
 & Thm.~\ref{thm:relative-virasoro-main}\\[4pt]

\multicolumn{3}{@{}l@{}}{\textbf{Compatibility with degeneration}}\\*[2pt]
$\mathcal X\to B$, $X_0$
 & Semistable family over a curve germ with smooth fiber $X$ and
   central fiber $X_0=X_1\cup_D X_2$.
 & Sec.~\ref{sec:degeneration}\\
$i_l$, $\iota_l$, $\iota_X$
 & Inclusions of $D$ into $X_l$, of $X_l$ into $\mathcal X$,
   and of the smooth fiber $X$ into $\mathcal X$.
 & \eqref{eq:general-degeneration}, \eqref{eq:matching-cohomology}\\
$H_\oplus$, $r$, $J$
 & The sum $H^*(X_1)\oplus H^*(X_2)$, the difference of restrictions
   to $D$, and the map $(i_{1*},-i_{2*})$ from $H^*(D)$.
 & \eqref{eq:matching-cohomology}\\
$K$, $\mathcal J$, $H$, $q$
 & Matching subspace $\ker r$, Gysin image $\operatorname{im}J$,
   quotient $H=K/\mathcal J$, and quotient map $q$.
 & \eqref{eq:matching-cohomology}\\
$\Psi$
 & Identification of the quotient $H$ with the smooth-fiber
   cohomology $H^*(X)$ under (C1).
 & Lem.~\ref{lem:complete-cohomology}\\
$\operatorname{Sp}_X$, $\operatorname{Sp}_{12}$
 & Novikov specializations from the smooth fiber or the two
   components to curve classes in $\mathcal X$.
 & \eqref{eq:general-degeneration}\\
$\Phi$, $\Phi_s$
 & Novikov specialization for deformation to the normal cone, and
   specialization for a complete degeneration using the lift $s$
   of the ordinary input and matching curve classes.
 & \eqref{eq:normal-cone-specialization}, \eqref{eq:general-specialization}\\

\end{longtable}
\endgroup

\bibliography{universal-BIB}
\bibliographystyle{amsxport}

\end{document}